\documentclass[opre,nonblindrev]{informs4aa} 

\OneAndAHalfSpacedXI 

\usepackage{tikz}
\usetikzlibrary{arrows.meta,decorations.pathreplacing}

\usepackage{natbib,booktabs}
 \bibpunct[, ]{(}{)}{,}{a}{,}{,}%
 \def\bibfont{\small}%
\usepackage{algorithm}
\usepackage{algpseudocode}
\usepackage{enumitem}
\usepackage{comment}
\usepackage{placeins}

\usepackage{xcolor}
\usepackage[normalem]{ulem}
\usepackage{hyperref}
\hypersetup{ 
 colorlinks,
 linkcolor={red!50!black},
 citecolor={blue!50!black},
 urlcolor={blue!80!black}
}
\usepackage{relsize}
\usepackage{bbm, cancel, cases}
\usepackage{subfigure}
\usepackage{algorithm}
\usepackage{algpseudocode}
\usepackage{tikz}
\usetikzlibrary{arrows.meta, decorations.pathreplacing, positioning, calc}
\usepackage{amsmath}
\TheoremsNumberedThrough
\ECRepeatTheorems

\renewcommand{\ORMSCLASSname}{%
    \textsf{\textbf{OR/MS subject classification:}}\enskip}

\usepackage{soul}
\setstcolor{gray}

\let\oldmarginpar\marginpar
\renewcommand\marginpar[1]{\-\oldmarginpar[\raggedleft\scriptsize \textcolor{red}{#1}]%
	{\raggedright\scriptsize \textcolor{red}{#1}}} 

\begin{document}

\RUNAUTHOR{A. Arlotto, Y. Liang, M. Liu, and Y. Wei}
\RUNTITLE{Online selection with disappearing types}

\TITLE{Online Selection with Disappearing Types: Protection Levels and Constant Optimality Gaps}

\ARTICLEAUTHORS{%
	\AUTHOR{Alessandro Arlotto}
	    \AFF{The Fuqua School of Business, Duke University, 100 Fuqua Drive, Durham, NC, 27708, \\ 
		\EMAIL{aa249@duke.edu}
        } 
	\AUTHOR{Yong Liang}
	    \AFF{Department of Management Science and Engineering, School of Economics and Management, Tsinghua University, Beijing 100084, China,\\  	\EMAIL{liangyong@sem.tsinghua.edu.cn}
        }	
    \AUTHOR{Min Liu}
	    \AFF{Department of Industrial Systems Engineering and Management, National University of Singapore, Singapore 117576, \\  			   
		\EMAIL{min.liu@nus.edu.sg}
        }
    \AUTHOR{Yehua Wei}
	    \AFF{The Fuqua School of Business, Duke University, 100 Fuqua Drive, Durham, NC, 27708, \\ 
		\EMAIL{yehua.wei@duke.edu}
        } 
} 

\ABSTRACT{
Motivated by arrival patterns studied in the protection-level literature in revenue management,
we consider online selection with demand types that are available only
during known, type-specific intervals of the planning horizon.
We show that the expected reward gap between the classical hindsight
offline benchmark and the optimal online policy can grow without bound
with the horizon.
We therefore introduce a tighter phase-aware benchmark that allows
hindsight within phases ending at type expirations while preserving
uncertainty about later phases.
We use this benchmark to construct two simple policies that combine
protection levels with budget-ratio control.
The phase protection policy uses a separate protection level for each
active phase--type pair.
The expiration protection policy uses only one protection level per type,
thereby reducing the number of protection levels that must be computed.
Under independent arrivals and a uniform positive lower bound on arrival
probabilities within active intervals, both policies achieve constant
optimality gaps, uniformly in the horizon and initial capacity for
fixed rewards and a fixed number of types.
}
\ORMSCLASS{
Primary: Dynamic programming/optimal control: Markov; Decision analysis: sequential.
Secondary: Inventory/production: policies; Analysis of algorithms: suboptimal algorithms.
}
\HISTORY{
This version: November 30, 2026.
}

\maketitle

\section{Introduction}

Online selection in capacity-control revenue management provides a framework for allocating scarce, often perishable capacity among sequentially arriving requests with heterogeneous values. Under this framework, firms must balance the immediate revenue from accepting a request against the opportunity cost of reducing the capacity available for uncertain future demand \citep{Talluri-VanRyzin-2006}.
Among the research streams informing these decisions, two are particularly relevant. First, a classical stream develops protection-level controls for allocating capacity across customer classes with different revenues and arrival patterns. 
Second, a more recent stream studies these decisions within the broader framework of stochastic online resource allocation, developing adaptive policies with constant-regret guarantees for (mostly) stationary arrivals. 

Protection-level controls reserve capacity for future demand by accepting a request only when the remaining capacity exceeds a specified threshold. 
In the classical two-fare model of \citet{Littlewood-1972}, low-fare demand arrives first and disappears before high-fare demand begins, and the protection level balances the current low fare against the expected marginal revenue from retaining capacity for later high-fare customers.
Subsequent work extended this approach to multiple fare classes \citep{Brumelle-McGill-1993}, settings without the low-before-high arrival assumption \citep{Lee-Hersh-1993,Robinson-1995}, and settings in which protection levels are learned from demand data \citep{Kunnumkal-Topaloglu-2009}.

The stochastic online resource-allocation literature develops another approach to capacity control. It focuses on bounding the regret, defined as the expected difference in total reward between the hindsight offline optimum and an online policy. 
Across a range of allocation settings with finitely many demand types and different arrival processes, budget-ratio and related adaptive policies achieve regret (hence, optimality gap) bounded uniformly in the time horizon and initial resource capacities \citep{Arlotto-Gurvich-2019,Vera-Banerjee-2020,vera2021online}.
However, the arrival patterns motivating classical protection-level controls, such as the example by \cite{Littlewood-1972},
pose a fundamental challenge for these constant-regret guarantees. In these settings, the information advantage of hindsight can make the offline benchmark too strong. As a result, even the optimal online policy may incur regret that grows with the horizon.

In this paper, we study online selection with demand types that are available only during known, type-specific intervals and disappear when those intervals end, an arrival pattern familiar from classical protection-level models.

Our contributions are threefold.

\begin{enumerate}[label=\textnormal{(\roman*)}]

\item \textbf{A phase-aware benchmark.}
In our setting, the standard hindsight offline benchmark can be too loose for establishing a constant optimality gap: its information advantage can create a gap from the optimal dynamic value that grows with the horizon 
(see Example~\ref{example:sqrt_T} for an instance with a horizon of $T$ periods in which even the optimal online policy incurs $\Omega(\sqrt{T})$ regret).
Hence, we introduce a tighter phase-aware information-relaxation benchmark that partitions the horizon into consecutive phases ending at the distinct type-expiration times. 
The benchmark is the optimal expected total reward under an information structure that allows hindsight within each phase but requires the capacity carried into future phases to be chosen before their arrivals are observed.

\item \textbf{Budget-ratio protection policies with $O(1)$ optimality gaps.}
We develop \emph{budget-ratio protection} policies that combine protection levels with budget-ratio control through a choice of future reference times. 
A protection level accounts for the value of capacity retained for demand after the reference time, while a certainty-equivalent demand adjustment based on expected arrivals guides allocation through that time. The resulting policy accepts a request when the remaining capacity is at least as large as the sum of these two quantities.
We develop two policies within this framework. 
The \emph{budget-ratio phase protection (PP)} policy uses the end of the current phase as its reference time and assigns a separate protection level to each active phase-type pair. 
By contrast, the \emph{budget-ratio expiration protection (EP)} policy uses each type's own expiration time throughout its active interval, requiring only one protection level per type.
The PP policy follows the benchmark’s protection structure, while the EP has a simpler parameterization with fewer protection levels but need not preserve that structure.
We show that both policies’ expected total rewards lie within $O(1)$ of the phase-aware benchmark, thereby establishing their constant optimality gaps.
These bounds are uniform over the horizon $T$ and initial capacity $B$, for a fixed number of types, fixed rewards, and a fixed positive lower bound on arrival probabilities within active intervals.

\item \textbf{A phase-by-phase comparison for the EP policy.}
Because the EP policy uses only one protection level per type, it need not match the benchmark’s phase-specific protection levels or preserve its reward-priority ordering (as the PP policy does).
To establish its constant optimality gap, we develop a phase-by-phase comparison that bounds the expected loss relative to the benchmark, accounting jointly for earned rewards and the benchmark continuation value of remaining capacity.
This comparison yields an upper bound with two components, one associated with arrival uncertainty and the other with differences between EP’s effective end-of-phase protection levels and the benchmark’s protection levels. 
We bound the first component using concentration inequalities. For the second component, 
we bound it through the protection-level value losses.
These losses are computed iteratively through phases from each type's expiration, allowing us to apply a backward propagation argument to obtain bounds uniform in the horizon and initial capacity.

\end{enumerate}

The remainder of the paper is organized as follows.
Section~\ref{sec:literature} reviews the related literature.
Section~\ref{sec:model} presents the model, introduces the phase-aware benchmark, and compares its value with the offline benchmark and the optimal dynamic value.
Section~\ref{sec:policy_guarantees} introduces the PP and EP policies and states their performance guarantees.
The main arguments supporting these guarantees appear in Section~\ref{sec:policy_analysis}, with technical details in the appendix.
Section~\ref{sec:numerical_block_noise} reports numerical experiments, and Section~\ref{sec:conclusions} offers concluding remarks. 

For any $z\in\mathbb R$, write $z^+:=\max\{z,0\}$. 
For function $V: \mathbb{R}^n \rightarrow \mathbb{R}$, we write $V\equiv0$ if $V$ is zero throughout its domain. We use $C_0$ to denote a finite constant independent of the horizon length $T$ and the initial capacity $B$.
For any positive integer $N$, we let $[N] := \{1, 2, \ldots, N\}$ denote the set of integers from $1$ through $N$.

\section{Related Literature}\label{sec:literature}

Our work is related to two streams of literature: online resource allocation with uniformly bounded loss and protection-level control in revenue management.

The first closely related stream studies online resource-allocation policies whose expected loss remains bounded uniformly in the horizon and the initial resource capacities. 
For the finite-support multisecretary problem with i.i.d. values, in which at most $B$ of $T$ sequentially observed values may be selected, \citet{Arlotto-Gurvich-2019} construct an adaptive budget-ratio policy with constant regret relative to the offline benchmark. \citet{Vera-Banerjee-2020} introduce compensated coupling and the Bayes Selector, obtaining constant-regret guarantees relative to suitable clairvoyant ex post benchmarks for a broader class of finite-type online-allocation problems. Follow-up work develops Bellman-inequality and model-predictive-control frameworks and studies the geometry and robustness underlying constant-regret guarantees
\citep{vera2021online,Banerjee-Freund-2024,Vera-et-al-2024}. For the multisecretary problem with more general support, \citet{Besbes-et-al-2024} establish distribution-dependent regret lower bounds that grow polynomially with the horizon for certain reward distributions. They develop policies whose regret matches these lower bounds up to logarithmic factors.
This literature also considers operational settings involving overbooking, bin-packing, delayed decisions, and online fulfillment
\citep{Freund-Zhao-2022,Banerjee-Freund-2024,Xie-Ma-Xin-2025,Arlotto-Keskin-Wei-2026,hssaine2026online}, as well as the fair allocation of perishable resources
\citep{banerjee2026online}.

We focus on online selection with types that disappear at different known times.
Although existing constant-regret results cover certain time-varying arrival processes, such guarantees relative to the offline benchmark are not generally attainable in our model class. Specifically, we exhibit instances in which the gap between the offline benchmark and the optimal dynamic policy is $\Omega(\sqrt{T})$. As a result, we instead establish an $O(1)$ gap relative to a phase-aware information-relaxation benchmark that upper-bounds the optimal dynamic value. The analysis of the expiration protection policy uses phase-by-phase
comparisons and a backward propagation argument.

A second related stream studies protection-level and booking-limit policies, which are widely used for capacity-control revenue management.
Specifically, the single-leg revenue management problem is a single-resource online allocation problem in which each accepted request consumes one unit of capacity and rewards take values from a finite set of fare classes.
Building on Littlewood's two-fare protection rule \citep{Littlewood-1972}, \citet{Brumelle-McGill-1993} characterize optimal fixed booking limits for multiple nested fare classes under low-before-high arrivals and independent fare-class demands. Dynamic formulations subsequently remove the fare-ordering assumption and characterize optimal policies through critical capacity values or booking limits \citep{Lee-Hersh-1993,Lautenbacher-Stidham-1999}. We refer the reader to \cite{Talluri-VanRyzin-2006} for a comprehensive review of the revenue management literature.

Other work develops protection-level and booking-limit policies under limited demand information, using historical observations or demand samples, bounds on fare-class demands, or predictive advice \citep{vanRyzin-McGill-2000,Kunnumkal-Topaloglu-2009,Lan-et-al-2008,Balseiro-et-al-2026,Golrezaei-et-al-2026}. Protection-level controls have also been extended to multiple-resource network revenue management through virtual nesting \citep{Bertsimas-deBoer-2005,vanRyzin-Vulcano-2008b}.

Our focus is on simple protection-level policies for online selection with type-specific active intervals.
We combine protection levels at selected future reference times with budget-ratio control based on expected arrivals through those times.
In particular, the expiration protection policy uses one protection level per type, with that type's expiration as its reference time.
We establish an $O(1)$ optimality gap relative to the optimal dynamic policy, uniformly in the horizon and initial capacity under our model assumptions.

\section{Model and Benchmark}\label{sec:model}

We consider an online allocation problem over a finite horizon of $T$ periods, indexed by $t \in [T]:=\{1,\ldots,T\}$.
There are $J$ types of arrivals, indexed by $j \in [J]$.
Each type $j$ is associated with a fixed reward $r_j>0$.
Without loss of generality, we label the types in (weakly) decreasing reward order, so that 
$r_1\ge r_2\ge\cdots\ge r_J$, and let $r_{\max}:=r_1$.
In each period $t \in [T]$, at most one arrival occurs.
Let $d_j(t) \in \{ 0, 1 \}$ denote the indicator that a type-$j$ arrival occurs in period $t$, so that $\sum_{j=1}^J d_j(t) \le 1$ for every $t \in [T]$.
We let $\lambda_{j,t} := \mathbb{E}[d_j(t)]$ be the probability of a type-$j$ arrival in period $t$, and note that the probability of no arrival in period $t$ is $1 - \sum_{j=1}^J \lambda_{j,t}$.
We assume that the arrival probabilities are known to the decision maker and that the arrival vectors $(d_1(t), \dots, d_J(t))$ are independent across periods. The decision maker has an initial capacity of $B \in \mathbb Z_+$ units.

In each period, if an arrival occurs, the decision maker first observes its type and then decides whether to accept or reject it before the next period begins.
Accepting a type-$j$ arrival consumes one unit of capacity and earns reward $r_j$.
Rejecting an arrival earns nothing and consumes no capacity.
All decisions are irrevocable and must be based only on information available at the time they are made.

An online policy $\pi$ specifies these decisions as functions of the known arrival probabilities, the current arrival outcome, and the history of past arrivals and actions.
Set $B_1^\pi = B$, and let $B_t^\pi$ denote the remaining capacity at the beginning of period $t$ under policy $\pi$.
We let $x_{j,t}^\pi = 1$ if policy $\pi$ accepts a type-$j$ arrival in period $t$, and $x_{j,t}^\pi = 0$ otherwise.
In each period, a feasible policy accepts only realized arrivals, and the number of accepted arrivals cannot exceed the remaining capacity:
\begin{align*}
    x_{j,t}^\pi \le d_j(t)
    \text{ for all } j \in [J], \ t \in [T],
    \qquad \text{and} \qquad
    \sum_{j=1}^J x_{j,t}^\pi \le B_t^\pi
    \text{ for all } t \in [T].
\end{align*}
The remaining capacity process then evolves as
\begin{align*}
    B_{t+1}^\pi
    = B_t^\pi - \sum_{j=1}^J x_{j,t}^\pi,
    \quad \forall t \in [T].
\end{align*}
The total reward earned by policy $\pi$ is
\begin{align*}
    \mathrm{Reward}(\pi)
    := \sum_{t=1}^T \sum_{j=1}^J r_j x_{j,t}^\pi.
\end{align*}

A key feature of our model is that types may become available after the horizon begins and disappear before it ends.
For each type $j \in [J]$, we define its \emph{activation} time $s_j$ and \emph{expiration} time $e_j$ as
\begin{equation}
    s_j := \min\left\{ t \in [T] : \lambda_{j,t} > 0 \right\}
    \qquad \text{and} \qquad
    e_j := \max\left\{ t \in [T] : \lambda_{j,t} > 0 \right\}.
    \label{eq:activation_expiration_times}
\end{equation}
We call $[s_j, e_j]$ the \emph{active interval} of type $j$.
By definition, $\lambda_{j,t} = 0$ for $t < s_j$ or $t > e_j$.
Thus, $e_j$ is the last period in which type $j$ may arrive.
Without loss of generality, we also assume that $\max_{j \in [J]} e_j = T$, since any trailing periods in which no type may arrive can be removed.

We impose the following \emph{active-interval assumption}: there exists a constant $\underline{\lambda} > 0$, independent of $T$ and $B$, such that
\begin{equation}
    \label{eq:arrival_probability_lower_bound}
    \lambda_{j,t} \ge \underline{\lambda},
    \qquad
    \text{for all $j \in [J]$ and all periods $t \in [s_j, e_j]$.}
\end{equation}
Under this assumption, a type's arrival probability may vary across periods within its active interval but must remain uniformly bounded away from zero.
Remark~\ref{rem:active-interval-assumption-relaxed} discusses alternative arrival models when \eqref{eq:arrival_probability_lower_bound} does not hold.

While this problem can in principle be solved by dynamic programming, with optimal expected total reward $V^\star$, our goal is to design simple online policies with provable performance guarantees.
Because a direct analysis of their performance relative to $V^\star$ can be difficult, a common approach is to compare their expected total rewards with a suitably chosen upper bound on $V^\star$ that has a simpler structure.
A typical choice is to use a full-information (offline) upper bound and analyze the regret of online policies.
Under independent arrivals and the active-interval assumption in \eqref{eq:arrival_probability_lower_bound}, with $s_j = 1$ and $e_j = T$ for all $j \in [J]$, simple online policies can achieve constant regret, 
as follows from the analysis of \citet[Theorem~3]{Vera-Banerjee-2020}.
Such a guarantee, however, need not extend to settings with disappearing types: even the optimal online policy can have regret that grows with the horizon, making the offline benchmark too coarse for establishing an $O(1)$ optimality gap.
We therefore seek to bound the \emph{optimality gap} of an online policy $\pi$, defined as the expected difference between the total rewards of an optimal dynamic programming policy and policy $\pi$, i.e., $V^\star - \mathbb E[\mathrm{Reward}(\pi)]$, using a tighter upper bound on $V^\star$.
To construct this benchmark, we partition the horizon into phases determined by type expirations and allow full information within each phase while retaining uncertainty about arrivals in later phases.
We motivate and formalize this benchmark in the following subsections.

\subsection{Insufficiency of the Offline Benchmark}\label{sec:phase_aware_benchmark}

Here we define the offline (full-information) benchmark and show why it can be too loose for establishing constant optimality gaps when types disappear.
The information advantage from observing the entire arrival sequence before allocating capacity can make its expected reward exceed that of the optimal online policy by $\Omega(\sqrt{T})$.

Given a realized arrival-count vector $\boldsymbol{A} = (A_1, \dots, A_J)$ and a capacity level $b \in \{ 0, 1, \dots, B \}$, let $P(b, \boldsymbol{A})$ denote the offline optimal reward obtainable from capacity $b$ when these arrival counts are known before allocating the available capacity.
This value is given by the linear program
\begin{align}
    P(b, \boldsymbol{A}) := \max_{x \in \mathbb{R}_+^J} \quad
    & \sum_{j=1}^J r_j x_j
    \label{eq:offline_value_def} \\
    \text{s.t.} \quad
    & \sum_{j=1}^J x_j \le b, \notag \\
    & 0 \le x_j \le A_j, \qquad \forall j \in [J]. \notag
\end{align}
Here, $x_j$ represents the number of type-$j$ arrivals accepted in hindsight.
Because each arrival consumes one unit of capacity and both capacity and arrival counts are integer-valued, this linear program admits an integral optimal solution.
Equivalently, $P(b, \boldsymbol{A})$ is the reward obtained by accepting, in hindsight, up to $b$ arrivals with the highest rewards among those in $\boldsymbol{A}$.
We then define the offline benchmark as
\begin{align*}
    \mathcal V_{\mathrm{off}}
    := \mathbb E[P(B, \boldsymbol A )].
\end{align*}
This benchmark is the conventional point of comparison in many online allocation problems.
The regret of an online policy $\pi$ is the difference $\mathcal V_{\mathrm{off}} - \mathbb E[\mathrm{Reward}(\pi)]$.

The following two-type example follows the classical low-before-high ordering in revenue management: all low-reward demand is realized before any high-reward demand.
It illustrates why the offline benchmark can be too loose to establish an $O(1)$ optimality gap for an online policy.

\begin{example}\label{example:sqrt_T}
    Consider a two-type online selection problem with rewards $r_1 = 1$ and $r_2 = \frac{ 1 }{ 2 }$ over a horizon of length $T$.
    For notational simplicity, assume that $T$ is divisible by $8$, and let the initial capacity be $B = 3 T / 8$.
    During the first half of the horizon, only type-$2$ arrivals may occur, with $\lambda_{1,t} = 0$ and $\lambda_{2,t} = \frac{ 1 }{ 2 }$ for $1 \le t \le \frac{ T }{ 2 }$.
    During the second half, only type-$1$ arrivals may occur, with $\lambda_{1,t} = \frac{ 1 }{ 2 }$ and $\lambda_{2,t} = 0$ for $\frac{ T }{ 2 } < t \le T$.
\end{example}

The following proposition gives an $\Omega(\sqrt{T})$ lower bound on the regret of every feasible online policy.

\begin{proposition}\label{prop:sqrt_T}
    For the family of two-type instances in Example~\ref{example:sqrt_T}, there is a constant $M > 0$, independent of $T$, such that every feasible online policy $\pi$ satisfies
    \begin{align*}
        \mathcal V_{\mathrm{off}} - \mathbb{E}[\mathrm{Reward}(\pi)]
        \ge \frac{ \sqrt{T} }{ 4 \sqrt{6} } - M.
    \end{align*}
\end{proposition}

In particular, $\mathcal V_{\mathrm{off}} - V^\star = \Omega(\sqrt{T})$ as $T \to \infty$ along this family of instances.
Thus, no feasible online policy can achieve a uniform $O(1)$ regret bound over our model class.

The lower bound has a simple newsvendor interpretation.
By the end of period $T/2$, the online policy has already committed to the capacity carried into the second half, without observing any high-reward demand.
Reserving too much capacity can leave capacity unused after rejecting low-reward arrivals, whereas reserving too little can force the rejection of high-reward arrivals.
Even an optimal choice of the capacity carried into the second half leaves an expected absolute mismatch of order $\sqrt{T}$ with that demand.

A related $\Omega(\sqrt{T})$ regret lower bound appears in \citet[Online Appendix C, Proposition 7]{Banerjee-Freund-2024}, where it establishes the sharpness of their assumption on the arrival process. 
Here, Example~\ref{example:sqrt_T} shows that permanent type disappearance can make the offline benchmark too loose to establish constant optimality gaps, motivating the tighter benchmark introduced next.

\subsection{A Phase-Aware Benchmark}

\begin{figure}[tbp]
    \FIGURE
    {\resizebox{\width}{!}{\begin{tikzpicture}[
        x=.88cm, y=.88cm,
        font=\fontsize{10}{12}\selectfont,
        >={Latex[length=1.5mm,width=1mm]},
        every node/.style={inner sep=2pt},
        interval/.style={
            fill=black!11,
            draw=black!65,
            line width=.4pt
        },
        endpoint/.style={
            black!55,
            dash pattern=on 2pt off 2.5pt,
            line width=.4pt
        },
        endpoint label/.style={
            font=\fontsize{6}{7}\selectfont,
            fill=black!11,
            inner xsep=1pt,
            inner ysep=0pt,
            outer sep=0pt
        }
    ]
    \foreach \j/\s/\e/\y in {1/2/12/3.1,2/3/7/2.1,3/1/10/1.1}{
        \node[anchor=east] at (.70,\y) {Type~\j\ ($r_{\j}$)};
        \draw[black!22,line width=.3pt] (1,\y) -- (12,\y);
        \draw[interval] (\s,{\y-.18}) rectangle (\e,{\y+.18});  
        \draw[black!70,line width=.55pt]
            (\s,{\y-.30}) -- (\s,{\y+.30)});
        \draw[black!70,line width=.55pt]
            (\e,{\y-.30}) -- (\e,{\y+.30});
        \node[endpoint label,anchor=west,xshift=2pt]
            at (\s,\y) {$s_{\j} = \s$};
        \node[endpoint label,anchor=east,xshift=-2pt]
            at (\e,\y) {$e_{\j} = \e$};
    }
    \foreach \e in {7,10,12}
        \draw[endpoint] (\e,0) -- (\e,3.50);
    \draw[->,line width=.5pt] (1,.26) -- (12.6,.26);
    \foreach \u in {1,...,12}{
        \draw[line width=.4pt] (\u,.21) -- (\u,.31);
        \node[anchor=north,font=\small,fill=white]
            at (\u,.15) {\u};
    }
    \node[anchor=west] at (12.64,.26) {$t$};
    \foreach \l/\r/\k in {1/7/1,8/10/2,11/12/3}{
        \draw[
            black!55,
            line width=.4pt,
            decorate,
            decoration={brace,mirror,amplitude=4pt}
        ] (\l - 0.25,-.36) -- (\r + 0.25,-.36);
        \node[anchor=north,font=\small]
            at ({(\l+\r)/2},-.64) {Phase~\k};
    }
    \end{tikzpicture}
    }}    
    {Type availability and phases\label{fig:model_timeline}}
    {Shaded intervals show when each type $j$ has a positive arrival probability.
    The labels $s_j$ and $e_j$ identify its activation and expiration times.
    Gray dashed lines mark the distinct type-expiration times $\tau_1 = 7$, $\tau_2 = 10$, and $\tau_3 = 12$. These define three phases, comprising periods 1--7, 8--10, and 11--12, respectively.}
\end{figure}

In Example~\ref{example:sqrt_T}, the offline solution uses knowledge of the second-half (high-reward) demand to determine how many low-reward arrivals to accept in the first half.
This informational advantage makes the offline benchmark too loose to establish a constant optimality gap.
We therefore construct a new benchmark that removes this advantage by allowing hindsight only within each half of the horizon, while requiring the capacity carried into the second half to be chosen before second-half demand is observed.

We extend this idea to the general model by partitioning the horizon into phases.
At the beginning of each phase, the benchmark solution uses the aggregate arrival counts within that phase, while arrivals in later phases remain unobserved.
We define the \emph{phase-aware benchmark}, denoted by $\mathcal V_{\mathrm{PA}}$, as the maximum expected total reward under this information structure.

Let $K \le J$ be the number of distinct type-expiration times, and denote these times by $\tau_1 < \cdots < \tau_K = T$.
Set $\tau_0 := 0$, and let phase $k$ consist of the periods $\mathcal T_k := \{ \tau_{k-1}+1, \ldots, \tau_k \}$.
Figure~\ref{fig:model_timeline} illustrates the active intervals and the corresponding phases of a problem instance with three types. 
Under the active-interval assumption, once a type is available during a phase, it remains available through the end of that phase.
Our analysis uses this continued availability to bound the expected loss within each phase by a constant independent of $T$ and $B$.
Since there are at most $J$ phases, summing these losses yields a bound that is also independent of $T$ and $B$.

The phase-aware benchmark can be computed by backward recursion over phases.
For each phase $k \in [K]$, let $A_j^{(k)} := \sum_{t \in \mathcal T_k} d_j(t)$ denote the number of type-$j$ arrivals during that phase, and let $\boldsymbol A^{(k)} := (A_1^{(k)}, \dots, A_J^{(k)})$ denote the corresponding arrival-count vector.
Let $\Psi_k(b)$ denote the optimal expected total reward from phase $k$ onward under the phase-aware information structure, with $b$ units of capacity at the beginning of phase $k$, before $\boldsymbol A^{(k)}$ is revealed.
Set $\Psi_{K+1} \equiv 0$.
Then, for $k \in [K]$ and $b \in \{ 0, \ldots, B \}$,
\begin{equation}
    \Psi_k(b)
    = \mathbb E_{\boldsymbol A^{(k)}}\left[
        \max_{\tilde b \in \{ 0, \ldots, b \}}
        \left\{ P(b - \tilde b, \boldsymbol A^{(k)}) + \Psi_{k+1}(\tilde b) \right\}
    \right],
    \label{eq:phase_aware_recursion}
\end{equation}
and the full-horizon benchmark value is $\mathcal V_{\mathrm{PA}} := \Psi_1(B)$.
After $\boldsymbol A^{(k)}$ is revealed, the inner maximization in \eqref{eq:phase_aware_recursion} determines how much capacity to allocate to the current phase and how much to retain for later phases.
Allocating $b - \tilde b$ units to the current phase gives the optimal reward $P(b - \tilde b, \boldsymbol A^{(k)})$, while retaining $\tilde b$ units gives the continuation value $\Psi_{k+1}(\tilde b)$.

The next lemma shows that the phase-aware benchmark upper-bounds $V^\star$ and, not surprisingly, is no larger than the offline benchmark.

\begin{lemma}[Bounding the phase-aware benchmark]
\label{lmm:benchmark_validity}
    The phase-aware benchmark satisfies
    \begin{align*}
        V^\star \le \mathcal V_{\mathrm{PA}} \le \mathcal V_{\mathrm{off}}.
    \end{align*}
\end{lemma}

Lemma~\ref{lmm:benchmark_validity} shows that the phase-aware benchmark can be used in lieu of the offline benchmark when establishing a constant optimality gap for a feasible online policy $\pi$.
Specifically, if $\pi$ satisfies
\begin{align*}
    \mathcal V_{\mathrm{PA}} - \mathbb E[\mathrm{Reward}(\pi)] \le C_0
\end{align*}
for some finite constant $C_0$ independent of $T$ and $B$, then its optimality gap is also at most $C_0$.

To guide the design and analysis of our online policies, we next characterize an optimal allocation for the inner maximization in \eqref{eq:phase_aware_recursion} using type-specific protection levels.

\begin{proposition}[Benchmark optimizer via protection levels]
\label{prop:direct_inner_optimizer}
    Fix a phase $k\in[K]$, and let $\mathcal J_k:=\{j\in[J]:s_j\le\tau_k\le e_j\}$ denote the types that can arrive during this phase. For each $j\in\mathcal J_k$, define the protection level
    \begin{equation}
        \theta_j(\tau_k)
        :=\max\left(\{0\}\cup\left\{u\in\{1,\ldots,B\}:\Psi_{k+1}(u)-\Psi_{k+1}(u-1)\ge r_j\right\}\right).
        \label{eq:theta_kj_def}
    \end{equation}
    Given capacity $b\in\{0,\ldots,B\}$ and a realized arrival vector $\boldsymbol A^{(k)}$, define
    \begin{align}
        x_j^*
        :=\min\left\{
            A_j^{(k)},
            \left(b-\theta_j(\tau_k)-\sum_{\substack{i\in\mathcal J_k\\ i<j}}A_i^{(k)}\right)^+
        \right\},\qquad j\in\mathcal J_k.
        \label{eq:xstar_direct}
    \end{align}
    Then $b^*:=b-\sum_{j\in\mathcal J_k}x_j^*$ is an optimal carry-over capacity,
    i.e., 
    $$b^* \in \arg\max_{\tilde b\in\{0,\ldots,b\}}\{P(b-\tilde b,\boldsymbol A^{(k)})+\Psi_{k+1}(\tilde b)\}.$$
    Moreover, $\{x_j^*\}_{j \in \mathcal J_k}$ is an optimal allocation for the inner problem $P(b-b^*,\boldsymbol A^{(k)})$.
\end{proposition}

The protection levels in Proposition~\ref{prop:direct_inner_optimizer} can be computed before the phase-$k$ arrival counts $\boldsymbol A^{(k)}$ are revealed.
The benchmark then uses these counts to determine the allocation $\{ x_j^* \}_{j \in \mathcal J_k}$ and carry-over capacity $b^*$.
This suggests that online policies using the benchmark's protection levels to account for the value of capacity retained beyond selected reference times can achieve expected total rewards close to the benchmark.
In the next section, we develop this idea into a class of policies that selects a protection reference time for each type-$j$ arrival in period $t$.
The reference time is a phase endpoint at or after $t$, and the corresponding protection level accounts for the value of capacity retained after that time.

\section{Protection Policies and Their Guarantees}\label{sec:policy_guarantees}

In this section, we present the budget-ratio protection framework, which combines the protection levels from Proposition~\ref{prop:direct_inner_optimizer} with the budget-ratio structure of \citet{Arlotto-Gurvich-2019}. 
The framework is specified in Algorithm~\ref{Algo:GeneralPolicy}.
For a type-$j$ arrival in period
$t$, the policy uses the end of a phase at or after $t$ as its protection reference time $\zeta_j(t)$. The corresponding
protection level $\theta_j(\zeta_j(t))$ accounts for the value of retaining
capacity for arrivals after that time.

Budget-ratio control then guides acceptance over $[t,\zeta_j(t)]$ as follows. Let $B_t$ denote the remaining capacity immediately before the period-$t$ decision. The policy compares the capacity above the protection level, $B_t-\theta_j(\zeta_j(t))$, with a demand adjustment defined as the sum of
the expected number of arrivals of types $i<j$ and half the expected number of type-$j$ arrivals over the interval $[t,\zeta_j(t)]$.
Intuitively, the first term in the sum estimates the future arrivals that are preferred over type $j$, 
while the latter term provides a buffer against arrival fluctuations. The policy accepts the arrival if $B_t\ge 1$ and the capacity above the protection level is at least as large as the demand adjustment.

\begin{algorithm}
    \caption{General Budget-Ratio Protection Policy}
    \label{Algo:GeneralPolicy}
    \begin{algorithmic}[1]
    \State \textbf{Input:} Known arrival probabilities
    $\{\lambda_{j,t}:j\in[J],\,t\in[T]\}$,
    protection reference times $\{\zeta_j(t)\}$ satisfying
    $\zeta_j(t)\in\{\tau_k:k\in[K],\ t\le\tau_k\le e_j\}$
    for $j\in[J]$ and $t\in [s_j,e_j]$,
    and initial capacity $B_1=B$.
    \For{$t=1,\ldots,T$}
        \State Observe the period-$t$ arrival outcome and set
        $x_{i,t}\gets0$ for all $i\in[J]$.
        \If{the arrival is of type $j$, $B_t \geq 1$, and \begin{align}
            B_t
            \ge
            \theta_j(\zeta_j(t))
            +
            \sum_{s=t}^{\zeta_j(t)}
            \left(
                \sum_{i<j}\lambda_{i,s}
                +\frac12\lambda_{j,s}
            \right).
            \label{eq:general_acceptance}
        \end{align}
    }
            \State $x_{j,t}\gets1$.
        \EndIf
        \State $B_{t+1}\gets B_t-\sum_{i=1}^Jx_{i,t}$.
    \EndFor
    \end{algorithmic}
\end{algorithm}

Next, we state two specific policies that fall within the budget-ratio protection framework, where the policies differ in their choice of protection reference times.
First, we choose the last period of the current phase as the protection reference time for every active type. Formally, this means setting $\zeta_j(t)=\tau_k$ whenever $t\in\mathcal T_k$ and type $j$ is active. We call this choice the \emph{budget-ratio phase protection} (PP) policy, as the protection reference time and protection levels for all types are reset after each phase. Let $\mathrm{Reward}(\mathrm{PP})$ denote its total reward. The following theorem states its performance guarantee.

\begin{theorem}[Performance guarantee for the PP policy]
\label{thm:direct_policy_main}
    Suppose the arrival probabilities satisfy the active-interval assumption in \eqref{eq:arrival_probability_lower_bound}. The phase protection policy (PP) satisfies
    \begin{align}
        0\le
        \mathcal V_{\mathrm{PA}}
        -
        \mathbb E[\mathrm{Reward}(\mathrm{PP})]
        \le
        J r_{\max}
        \left(
            \left\lceil\frac{4}{\underline\lambda}\right\rceil
            +\frac{16}{\underline\lambda^2}
        \right).
        \label{eq:direct_main_bound}
    \end{align}
\end{theorem}

For fixed $J$, $r_{\max}$, and $\underline\lambda > 0$,
Theorem~\ref{thm:direct_policy_main} provides an $O(1)$ bound
uniform over the horizon $T$, the initial capacity $B$, the phase lengths,
and all arrival-probability sequences satisfying the active-interval assumption.
Since $V^\star \le \mathcal V_{\mathrm{PA}}$, it also establishes
an $O(1)$ optimality gap relative to the optimal dynamic policy.
The analysis of this result is presented in
Subsection~\ref{subsec:direct_policy_analysis}.

The PP policy requires a separate protection level for each
phase--type pair. With $j \in \mathcal J_k$ and $k \in [K]$,
it may require as many as $KJ$ protection parameters in total.
In particular, if a type remains active after a phase ends,
the policy requires a new protection level for that type in the next phase.
One might wonder whether a single protection level per type
can suffice to retain an $O(1)$ optimality gap, thereby reducing
the number of required protection levels to $J$.
The answer is affirmative: we design the
\emph{budget-ratio expiration protection} (EP) policy,
which achieves this reduction while retaining an $O(1)$ optimality gap.

Formally, the EP policy sets $\zeta_j(t)=e_j$ throughout type $j$'s active interval and assigns a single protection level to that type.
Let $\mathrm{Reward}(\mathrm{EP})$ denote its total reward.
The following theorem establishes its performance guarantee.

\begin{theorem}[Performance guarantee for the EP policy]
    \label{thm:heterogeneous_EA_performance}
    Suppose the arrival probabilities satisfy the active-interval assumption in \eqref{eq:arrival_probability_lower_bound}. Then the expiration protection (EP) policy satisfies
    \begin{align}
        \mathcal V_{\mathrm{PA}}-\mathbb E[\mathrm{Reward}(\mathrm{EP})]
        \le{}&
        K
        \left(
            1 + \frac{8}{3 \underline \lambda} + \frac{2}{\exp(\underline \lambda^2/32)-1}
        \right)r_{\max} +
        \frac{
            32r_{\max}JK(K-1)
        }{
            3\underline\lambda^3
        }.
        \label{eq:heterogeneous_explicit_bound}
    \end{align}
\end{theorem}

Since $K\le J$, Theorem~\ref{thm:heterogeneous_EA_performance} gives a bound independent of $T$, $B$, and the phase lengths for fixed $J$, $r_{\max}$, and $\underline\lambda$.
Because $V^\star\le\mathcal V_{\mathrm{PA}}$, it also establishes an $O(1)$ optimality gap.

Next, we briefly describe how we establish the performance guarantees
of the PP and EP policies, and why the analysis of the EP policy is significantly
more difficult. 
First, the PP policy uses the benchmark's end-of-phase protection
levels $\theta_j(\tau_k)$.
These levels account for the value of capacity saved for later phases.
Within the current phase, PP uses expected arrival counts to guide its
decisions, while the benchmark knows the realized counts.
The buffer in the demand adjustment protects against this uncertainty:
when sufficiently many periods remain, a PP decision can lose value
relative to the benchmark only if the remaining arrival counts are far
from their expectations.
Concentration inequalities give error probabilities that decrease
exponentially with the number of remaining periods.
Summing these probabilities over the phase gives a bound independent
of the phase length.
Since each decision's loss is at most $r_{\max}$, the expected loss
per phase is $O(1)$.

The analysis for the EP policy is more delicate, as there is an additional
source of loss through EP's effective end-of-phase protection levels,
which may differ from those of the benchmark.
To compare EP with PP within
phase $k$, we split EP's demand adjustment at the phase end. For each
type $j\in\mathcal J_k$, define the \emph{effective end-of-phase protection level}
\begin{equation}\label{eq:ep_eop_protection_level}
    \eta_j^{\tau_k}
    :=
    \theta_j(e_j)
    +
    \sum_{s=\tau_k+1}^{e_j}
    \left(
        \sum_{i<j}\lambda_{i,s}
        +\frac12\lambda_{j,s}
    \right),
\end{equation}
which combines the expiration-time protection level with the demand
adjustment for periods after phase $k$ through type $j$'s expiration. Within phase $k$, EP uses the same demand adjustment as PP, with
$\eta_j^{\tau_k}$ in place of $\theta_j(\tau_k)$ to protect capacity for
future phases. Their acceptance thresholds therefore differ by
$\eta_j^{\tau_k}-\theta_j(\tau_k)$. This deterministic difference can lead
to suboptimal decisions even when remaining arrivals are close to their
expectations. Therefore, controlling the resulting loss requires an
additional argument beyond concentration of the current phase's arrivals. 
In Subsection~\ref{subsec:heterogeneous_performance}, we compare how many selected rewards and continuation marginal values exceed each cutoff under the benchmark and the policy.
Integrating these count differences gives the value gap within a phase, and a backward-propagation argument controls the loss associated with these differences in protection levels.

\begin{remark}
\label{rem:active-interval-assumption-relaxed}
Theorems \ref{thm:direct_policy_main} and \ref{thm:heterogeneous_EA_performance} derive sub-optimality guarantees under the active-interval assumption specified as \eqref{eq:arrival_probability_lower_bound}. When such assumption does not hold, one way to derive sub-optimality guarantees is to approximate them by modified probabilities $\{\tilde\lambda_{j,t}:j\in[J],\,t\in[T]\}$. 
Specifically, $\tilde\lambda_{j,t}$ denotes the probability of a type-$j$ arrival in period $t$ in an approximating model that satisfies the assumption. We measure the difference between the original and approximating arrival models by the \emph{total approximation error}
\begin{align*}
    \varepsilon
    :=
    \sum_{t=1}^T\sum_{j=1}^J
    |\lambda_{j,t}-\tilde\lambda_{j,t}|.
\end{align*}
This quantity sums the absolute changes in arrival probabilities over all types and periods. Appendix~\ref{app:arrival_approximation} bounds the effect of this approximation on the expected rewards of our policies. Specifically, if either policy is constructed using the modified probabilities and then kept fixed, its expected total rewards under the two arrival models differ by at most $r_{\max}\varepsilon$. Thus, this difference is bounded independently of $T$ and $B$ when $\varepsilon$ is uniformly bounded. If $\varepsilon$ grows slowly as $T$ and $B$ increase, the bound would grow proportionally to $\varepsilon$.
\end{remark}

\subsection{Connection to Classical Protection Levels and Budget-Ratio Policies}

We now interpret our budget-ratio protection policies in terms of their connections to classical protection levels and budget-ratio policies.

The protection level $\theta_j(\zeta_j(t))$ reserves capacity for arrivals after the protection reference time $\zeta_j(t)$.
This is similar to the classical single-resource revenue management models, where protection levels are chosen by comparing the reward of the current class with the marginal value of retaining capacity for future demand. For example, in Littlewood's two-class model, class $2$ arrives before class $1$ and $r_1>r_2$.
Let $e_2$ denote the last period of class-2 demand and let $A_1$ denote the subsequent class-1 demand.
With protection reference time $e_2$, the corresponding protection level is
$$
    \theta_2(e_2)
    =
    \max\left(
        \{0\}
        \cup
        \left\{
            u\in\{1,\ldots,B\}:
            r_1\mathbb P(A_1\ge u)\ge r_2
        \right\}
    \right).
$$
This is Littlewood's protection level for class $1$, up to the discrete tie-breaking convention. 
In this example, our protection level coincides with Littlewood's because both use the same marginal-value comparison. However, unlike Littlewood's two-class model, our setting does not require low-reward demand to precede high-reward demand.

The demand adjustment in \eqref{eq:general_acceptance} guides how capacity above the protection level is used from the current period through the reference time.
This part of the acceptance rule follows the comparison structure of the budget-ratio policy of \citet{Arlotto-Gurvich-2019}, which compares remaining capacity per remaining period with demand thresholds determined by the type probabilities. Specifically, subtracting the protection level and dividing by the number of periods in $[t,\zeta_j(t)]$, we can rewrite our acceptance threshold comparison as
\begin{align*}
    \frac{B_t-\theta_j(\zeta_j(t))}{\zeta_j(t)-t+1}
    \ge
    \sum_{i<j}
    \frac{\sum_{s=t}^{\zeta_j(t)}\lambda_{i,s}}{\zeta_j(t)-t+1}
    +
    \frac12
    \frac{\sum_{s=t}^{\zeta_j(t)}\lambda_{j,s}}{\zeta_j(t)-t+1}.
\end{align*}
This comparison matches the budget-ratio policy defined in \citet{Arlotto-Gurvich-2019} in their stationary model, with $\zeta_j(t)=T$ and $\theta_j(T)=0$ for all types, since $T$ is at the end of the horizon.

\section{Analysis of the Protection-Level Policies}
\label{sec:policy_analysis}

In this section, we present the main arguments and technical results used to establish the performance guarantees in Theorems~\ref{thm:direct_policy_main} and~\ref{thm:heterogeneous_EA_performance}.
We bound the policies' expected reward gaps relative to the phase-aware benchmark.
For the PP policy, its use of the benchmark's phase-specific protection levels allows us to bound the expected loss by controlling the probability of a tracking error in each period.
For the EP policy, we use a phase-by-phase comparison that accounts jointly for earned rewards and the benchmark continuation value of remaining capacity, together with a backward-propagation argument.

\subsection{Analysis of the Phase Protection Policy}
\label{subsec:direct_policy_analysis}

In this subsection, we analyze the PP policy's expected reward gap
relative to the phase-aware benchmark.
Throughout, we use $B_t^{\mathrm{PP}}$ to denote its remaining capacity
at the beginning of period $t$.
We apply the compensated coupling argument of \citet{Vera-Banerjee-2020}
within each phase.
At each period, we compare the policy's decision with an optimal
decision starting from the same remaining capacity, using hindsight
about the remaining arrivals in the current phase and the benchmark
continuation value for capacity carried into later phases.
We call a policy decision that is suboptimal for this optimization
problem a \emph{tracking error}.
We bound the probabilities of these errors and sum their associated
losses to obtain a constant expected loss per phase.

Fix any $t\in\mathcal T_k$, let $\widehat A_j^{(k)}(t):=\sum_{s=t}^{\tau_k}d_j(s)$ denote the number of type-$j$ arrivals remaining in phase $k$, including period $t$, and let $\widehat{\boldsymbol A}^{(k)}(t) :=(\widehat A_j^{(k)}(t))_{j\in[J]}$. 
Then, following the inner maximization in the phase-aware benchmark recursion \eqref{eq:phase_aware_recursion}, we define the \emph{periodic} maximum phase-aware expected reward at period $t$ under capacity $b$ as
\begin{align}
    W_{k,t} \bigl(b,\widehat{\boldsymbol A}^{(k)}(t)\bigr)
    :=
    \max_{\tilde b\in\{0,\ldots,b\}}
    \left\{
        P\bigl(b-\tilde b,\widehat{\boldsymbol A}^{(k)}(t)\bigr)
        +
        \Psi_{k+1}(\tilde b)
    \right\},
    \label{eq:direct_W_EX}
\end{align}
with $W_{k,\tau_k+1}(b,\mathbf 0)=\Psi_{k+1}(b)$.
When $t = \tau_{k-1}+1$, given a realized sequence of arrivals and capacity $b$, the maximization problem in \eqref{eq:phase_aware_recursion} is exactly $W_{k,t}\bigl(b,\widehat{\boldsymbol A}^{(k)}(t)\bigr)$, implying that 
\begin{align*}
 \Psi_k(b)
    = \mathbb E_{\boldsymbol A^{(k)}} \left[W_{k,\tau_{k-1}+1}\bigl(b,\widehat{\boldsymbol A}^{(k)}(t)\bigr)\right]
\end{align*}
Moreover, conditioned on $d_j(t) =  1$, we have
\begin{align}
W_{k,t}\bigl(b,\widehat{\boldsymbol A}^{(k)}(t)\bigr) = \max_{\substack{a\in\{0,1\}\\a\le b}}
    \left\{
        r_ja
        +
        W_{k,t+1}
        \bigl(b-a,\widehat{\boldsymbol A}^{(k)}(t+1)\bigr)
    \right\}.
\label{eq:satisfying_equation}
\end{align}
Therefore, for a type-$j$ arrival in period $t$, we call a feasible accept/reject decision a \emph{satisfying action} if it is optimal for the optimization problem above, and use $\mathcal S_{j,t}^{(k)}(b)$ to denote the set of all satisfying actions.

For $t\in\mathcal T_k$, we say that the PP policy makes a \emph{tracking error} for type $j$ if a type-$j$ arrival occurs and the policy's decision is not optimal for \eqref{eq:satisfying_equation}.
The corresponding event is defined as
\begin{align}
    \mathrm{Error}_{j,t}^{(k)}
    :=
    \left\{
        d_j(t)=1,\;
        x_{j,t}^{\mathrm{PP}}
        \notin
        \mathcal S_{j,t}^{(k)}(B_t^{\mathrm{PP}})
    \right\},
    \label{eq:error_direct_phasek}
\end{align}
where $x_{j,t}^{\mathrm{PP}}$ denotes the policy's acceptance indicator. The following proposition bounds the probability of this event.

\begin{proposition}[Tracking error probability]
\label{prop:direct_tracking_error}
Fix phase $k$, period $t\in\mathcal T_k$, and type $j\in\mathcal J_k$ with $\lambda_{j,t}>0$. If $\tau_k-t+1\ge4/\underline\lambda$, then
\begin{align}
    \mathbb P
    \left(
        \mathrm{Error}_{j,t}^{(k)}
        \,\middle|\,
        d_j(t)=1
    \right)
    \le
    2\exp\left(
        -\frac{\underline\lambda^2}{8}(\tau_k-t+1)
    \right).
    \label{eq:direct_tracking_error_bound}
\end{align}
Moreover, the event of a tracking error in period $t$,
$\mathrm{Error}_t:=\bigcup_{j\in[J]}\mathrm{Error}_{j,t}^{(k)}$, satisfies
\begin{align}
    \mathbb P\left(\mathrm{Error}_t\right)
    \le
    2\exp\left(
        -\frac{\underline\lambda^2}{8}(\tau_k-t+1)
    \right).
    \label{eq:union_tracking_error_bound}
\end{align}
\end{proposition}

Proposition~\ref{prop:direct_tracking_error} shows that the PP policy selects an optimal action for the benchmark optimization problem with high probability when many periods remain in the phase.
We are now ready to prove Theorem~\ref{thm:direct_policy_main}.

\proof{Proof of Theorem~\ref{thm:direct_policy_main}.}
For the sake of notational simplicity, write $W_{k,t}(b)$ for
$W_{k,t}(b,\widehat{\boldsymbol A}^{(k)}(t))$, the \emph{periodic} maximum phase-aware expected reward at period $t$ under capacity $b$.
For each $t\in\mathcal T_k$, define the one-period loss from the PP policy's decision by
\begin{align*}
    \ell_t
    :=
    W_{k,t}(B_t^{\mathrm{PP}})
    -\sum_{j=1}^J r_j x_{j,t}^{\mathrm{PP}}
    -W_{k,t+1}(B_{t+1}^{\mathrm{PP}}).
\end{align*}
If no arrival occurs, $\ell_t=0$. When an arrival occurs, \eqref{eq:satisfying_equation} gives $\ell_t$ as the gap between the optimal action value and the policy's action value. This gap is nonnegative and satisfies $\ell_t \in [0,r_{\max}]$, since removing one unit costs at most one arrival reward or one continuation marginal value
(Lemma~\ref{lem:psi_discrete_concavity}). Moreover, $\ell_t>0$ exactly when a tracking error occurs in period $t$.

Summing $\ell_t$ over all $t$ within phase $k$ cancels all of the intermediate $W_{k,t}$ terms except $t =\tau_{k-1}+1$ and $t = \tau_{k}+1$.
Since the phase-$k$ arrivals are independent of $B_{\tau_{k-1}+1}^{\mathrm{PP}}$,
the conditional expectation of
$W_{k,\tau_{k-1}+1}(B_{\tau_{k-1}+1}^{\mathrm{PP}})$ (over phase-$k$ arrivals)
is $\Psi_k(B_{\tau_{k-1}+1}^{\mathrm{PP}})$ by
\eqref{eq:phase_aware_recursion}. Thus,
\begin{align*}
    \mathbb E\!\left[
        \Psi_k(B_{\tau_{k-1}+1}^{\mathrm{PP}})
        - \sum_{t\in\mathcal T_k}\sum_{j=1}^J
r_jx_{j,t}^{\mathrm{PP}}
        -\Psi_{k+1}(B_{\tau_k+1}^{\mathrm{PP}})
    \right] =\sum_{t\in\mathcal T_k}\mathbb E[\ell_t]
\end{align*}

Summing this equality over all $K$ phases cancels the intermediate continuation-value terms.
Because $\Psi_1(B)=\mathcal V_{\mathrm{PA}}$ and
$\Psi_{K+1}\equiv0$, we have
\begin{equation}\label{eq:telescope}
\mathcal V_{\mathrm{PA}}-\mathbb E[\mathrm{Reward}(\mathrm{PP})] = \sum_{k \in [K]} \sum_{t\in \mathcal T_k}\mathbb E[\ell_t]. 
\end{equation}

Since a tracking error causes a loss of at most $r_{\max}$, for each phase $k$,
\begin{equation}    
\label{eq:loss_by_phase}
    \sum_{t\in\mathcal T_k}\mathbb E[\ell_t]
    \le r_{\max}\sum_{t\in\mathcal T_k}\mathbb P(\mathrm{Error}_t) 
    \le r_{\max}
    \left(\left\lceil\frac4{\underline\lambda}\right\rceil
        +2\sum_{h=\lceil4/\underline\lambda\rceil+1}^{\infty}
        e^{-\underline\lambda^2h/8}\right),
\end{equation}
where the last inequality follows from Proposition~\ref{prop:direct_tracking_error}. More specifically, for $h=\tau_k-t+1$, we bound the error probability by one if
$h\le\lceil4/\underline\lambda\rceil$ and apply
\eqref{eq:union_tracking_error_bound} otherwise.
Finally, combining \eqref{eq:telescope} and \eqref{eq:loss_by_phase} yields
\begin{align*}
\mathcal V_{\mathrm{PA}}-\mathbb E[\mathrm{Reward}(\mathrm{PP})]
\le \sum_{k \in [K]}\sum_{t \in \mathcal T_k} \mathbb E[\ell_t] 
        &\le K r_{\max}
    \left(
        \left\lceil\frac4{\underline\lambda}\right\rceil
        +2\sum_{h=\lceil4/\underline\lambda\rceil+1}^{\infty}
        e^{-\underline\lambda^2h/8}
    \right)\\
    &\le J r_{\max}
    \left(
        \left\lceil\frac4{\underline\lambda}\right\rceil
        +\frac{16}{\underline\lambda^2}
    \right),
\end{align*}
where the last inequality uses $K\le J$ and $\sum_{h=1}^{\infty}e^{-xh}\le1/x$ for $x>0$. 
\Halmos\endproof

\subsection{Analysis of the Expiration Protection Policy}
\label{subsec:heterogeneous_performance}

In this subsection, we present the main steps for establishing Theorem~\ref{thm:heterogeneous_EA_performance} through a phase-by-phase comparison of the EP policy with the phase-aware benchmark.
As discussed in Section~\ref{sec:policy_guarantees}, the EP policy's effective end-of-phase protection levels may differ from the benchmark's phase-specific protection levels.
These differences can cause suboptimal decisions even when arrivals are close to their expectations, so the concentration argument for the PP policy 
to bound the tracking error probability in each period
is not sufficient on its own to establish this guarantee. However, a large difference in protection levels need not imply a large loss in total reward. If a type remains active in a later phase, capacity saved by rejecting an arrival may allow an additional acceptance of that type, offsetting the shortfall in current-phase reward. We therefore develop a phase-by-phase comparison that jointly accounts for earned rewards and the benchmark continuation value of remaining capacity. 

Our analysis of the EP policy is divided into the following steps.
First, in Section~\ref{subsubsec:heterogeneous_value_representation}, we use a layer-cake representation to express the value gap within a phase between the benchmark and an arbitrary feasible policy as an integral of differences in the numbers of selected arrival rewards and continuation marginal values above each cutoff.
Second, in Section~\ref{subsubsec:heterogeneous_pointwise_bound}, we follow this representation and bound the expected difference between the two counts at each cutoff, accounting separately for arrival uncertainty and differences between the effective end-of-phase protection levels and the benchmark's phase-specific protection levels.
Finally, in Section~\ref{subsubsec:heterogeneous_backward_propagation}, we integrate this bound over cutoffs and bound the contribution from protection-level differences by a sum of type-specific value losses. We control each type's loss by backward propagation from its expiration, where the loss is zero. Summing the resulting bounds completes the proof of Theorem~\ref{thm:heterogeneous_EA_performance}.

\subsubsection{Representation of the One-Phase Value Gap}
\label{subsubsec:heterogeneous_value_representation}

In this subsection, we use a layer-cake representation to express the value gap within a phase between the benchmark and an arbitrary feasible policy as an integral of differences in the numbers of selected arrival rewards and continuation marginal values above each cutoff.

Fix a phase $k$, an entering capacity $b\in\{0,\ldots,B\}$, and an arbitrary feasible policy $\pi$ during this phase. Let $R_k^\pi$ denote its reward during the phase and $b'$ its remaining capacity at the end of the phase. We evaluate the policy's allocation by adding the reward earned during the phase to the benchmark continuation value of its remaining capacity, giving $R_k^\pi+\Psi_{k+1}(b')$.

The reward $R_k^\pi$ is a sum of rewards from accepted arrivals. To express the continuation value in the same form, we decompose it into marginal values, one for each unit of capacity carried forward as
\begin{align}
    \nu_u^{(k)}
    :=
    \Psi_{k+1}(u)-\Psi_{k+1}(u-1),
    \qquad u=1,\ldots,B.
    \label{eq:heterogeneous_continuation_marginals}
\end{align}
These marginal values satisfy $r_{\max}\ge\nu_1^{(k)}\ge\cdots\ge\nu_B^{(k)}\ge0$. Since $\Psi_{k+1}(0)=0$, we have
\begin{align*}
    R_k^\pi+\Psi_{k+1}(b')
    =
    \sum_{t\in\mathcal T_k}\sum_{j\in[J]}
    r_jx_{j,t}^{\pi}
    +
    \sum_{u=1}^{b'}\nu_u^{(k)},
\end{align*}
where $x_{j,t}^{\pi}$ denotes the type-$j$ acceptance indicator
under policy $\pi$ in period $t$. Applying the so-called layer-cake representation to this sum gives
\begin{align}
    R_k^\pi+\Psi_{k+1}(b')
    =
    \int_0^{r_{\max}}N^{\pi}(v)\,dv,
    \label{eq:heterogeneous_policy_representation}
\end{align}
where, for each cutoff $v\in[0,r_{\max}]$,
\begin{align}
    N^{\pi}(v)
    :=
    \sum_{t\in\mathcal T_k}
    \sum_{\substack{j\in[J]\\r_j>v}}
    x_{j,t}^{\pi}
    +
    \sum_{u=1}^{b'}
    \mathbf 1\{\nu_u^{(k)}>v\}.
    \label{eq:heterogeneous_Npi_definition}
\end{align}
Here, $N^{\pi}(v)$ counts the accepted arrival rewards and continuation marginal values that are strictly greater than $v$.

Similarly, given the realized arrivals in phase $k$, the benchmark maximizes the same sum of the reward earned during the phase and the continuation value of the capacity carried forward over feasible allocations of the $b$ units of capacity. Since the continuation marginal values are nonincreasing, this maximization is equivalent to selecting the $b$ largest elements from the multiset containing the rewards of all arrivals realized during phase $k$ and the continuation marginal values $\nu_1^{(k)},\ldots,\nu_b^{(k)}$. Applying the layer-cake representation to this selection gives
\begin{align}
    \max_{\tilde b\in\{0,\ldots,b\}}
    \left\{
        P(b-\tilde b,\boldsymbol A^{(k)})
        +\Psi_{k+1}(\tilde b)
    \right\}
    =
    \int_0^{r_{\max}}N^*(v)\,dv,
    \label{eq:heterogeneous_benchmark_representation}
\end{align}
where
\begin{align}
    N^*(v)
    :=
    \min\left\{
        b,\,
        \sum_{\substack{j\in[J]\\r_j>v}}A_j^{(k)}
        +
        \sum_{u=1}^{b}\mathbf 1\{\nu_u^{(k)}>v\}
    \right\}
    \label{eq:heterogeneous_Nstar_definition}
\end{align}
counts the elements above $v$ in the benchmark's selection.

Combining the two layer-cake representations \eqref{eq:heterogeneous_policy_representation} and \eqref{eq:heterogeneous_benchmark_representation}, we have
\begin{align}
        &\max_{\tilde b\in\{0,\ldots,b\}}
        \left\{
            P(b-\tilde b,\boldsymbol A^{(k)})
            +\Psi_{k+1}(\tilde b)
        \right\}
        -R_k^\pi-\Psi_{k+1}(b')
        =
        \int_0^{r_{\max}}
        \left[N^*(v)-N^{\pi}(v)\right]\,dv,
        \label{eq:heterogeneous_item_representation}
    \end{align}
which expresses the one-phase value gap by comparing the numbers of elements greater than each cutoff in the benchmark's and the policy's selections.
We do not need to identify, for each element selected by the benchmark but omitted by the policy, which element the policy selects in its place. A smaller number of accepted arrivals with rewards above a cutoff can be offset by a larger number of selected continuation marginal values above that cutoff.
The representation therefore accounts for their combined contribution, which may be small 
even when the policy's acceptance decisions are significantly different from the benchmark.

\subsubsection{Separating Arrival Uncertainty and Protection-Level Differences}
\label{subsubsec:heterogeneous_pointwise_bound}

In this subsection, we bound the expected difference $\mathbb E[N^*(v)-N^{\mathrm{EP}}(v)]$ between the benchmark and policy counts at each cutoff, separating the effects of arrival uncertainty and differences in protection levels. We take expectations over the arrivals during phase $k$, with the entering capacity $b$ fixed.

We first measure the effect of using the effective end-of-phase protection levels $\eta_j^{\tau_k}$ (defined in \eqref{eq:ep_eop_protection_level})
instead of the benchmark's protection levels $\theta_j(\tau_k)$. At a cutoff $v\in[0,r_{\max}]$, we compare each effective level, rounded up and capped at $b$, with the number of continuation marginal values above $v$ among $\nu_1^{(k)},\ldots,\nu_b^{(k)}$. For a type with $r_j>v$, protecting more units than this number may cause the rejection of an arrival with reward above $v$ while retaining units whose continuation marginal values are at most $v$. For a type with $r_j\le v$, protecting fewer units may allow an accepted arrival to displace a continuation marginal value above $v$. We measure the largest of these excesses or shortfalls by
\begin{align}
    \xi(v):=\max\Bigg\{
        &\max_{\substack{j\in\mathcal J_k\\r_j>v}}
        \left(
            \min\{\lceil\eta_j^{\tau_k}\rceil,b\}
            -\sum_{u=1}^{b}\mathbf 1\{\nu_u^{(k)}>v\}
        \right)^+,
        \notag\\
        &\max_{\substack{j\in\mathcal J_k\\r_j\le v}}
        \left(
            \sum_{u=1}^{b}\mathbf 1\{\nu_u^{(k)}>v\}
            -\min\{\lceil\eta_j^{\tau_k}\rceil,b\}
        \right)^+
    \Bigg\},
    \label{eq:heterogeneous_d_definition}
\end{align}
where a maximum over an empty set is zero. The quantity $\xi(v)$ also depends on the phase $k$ and entering capacity $b$; we omit these indices for brevity. With $k$ and $b$ fixed, $\xi(v)$ is a deterministic function of $v$.

The benchmark's phase-specific protection levels satisfy
\begin{align*}
    \min\{\theta_j(\tau_k),b\}
    =
    \sum_{u=1}^{b}\mathbf 1\{\nu_u^{(k)}\ge r_j\}.
\end{align*}
Consequently, replacing $\eta_j^{\tau_k}$ with $\theta_j(\tau_k)$ in \eqref{eq:heterogeneous_d_definition} gives $\xi(v)=0$ for every cutoff. This reflects the benchmark's marginal-value comparison: its phase-specific protection level for type $j$ protects exactly the units whose continuation marginal values are at least $r_j$.

To describe arrival uncertainty at cutoff $v$, define
\begin{equation}
    \Delta_v(t)
    :=\sum_{s=t+1}^{\tau_k}\sum_{\substack{j\in[J]\\r_j>v}}
    \bigl(d_j(s)-\lambda_{j,s}\bigr),
    \qquad t\in\mathcal T_k.
    \label{eq:heterogeneous_arrival_deviation}
\end{equation}
Thus, $\Delta_v(t)$ is the realized number of arrivals with rewards above $v$ from period $t+1$ through the end of phase $k$, minus its expectation. The following lemma bounds the difference between the two counts using $\xi(v)$ and these arrival deviations.

\begin{lemma}[Separating arrival uncertainty and protection-level differences]
    \label{lem:heterogeneous_pointwise_count}
    Suppose the arrival probabilities satisfy the active-interval assumption in \eqref{eq:arrival_probability_lower_bound}. 
    Fix a phase $k$ and an entering capacity $b\in\{0,\ldots,B\}$, and run the EP policy during phase $k$ from capacity $b$. For every cutoff $v\in[0,r_{\max}]$ and every realization of the phase arrivals,
    \begin{align}
        N^*(v)-N^{\mathrm{EP}}(v)
        \le{}&
        \left\lceil\frac{8(\xi(v)+1)}{3\underline\lambda}\right\rceil
        +\sum_{t\in\mathcal T_k}
        \mathbf 1\left\{
            |\Delta_v(t)|>\frac{\underline\lambda}{8}(\tau_k-t)
        \right\}.
        \label{eq:heterogeneous_pointwise_separation}
    \end{align}
    Consequently,
    \begin{align}
        \mathbb E\!\left[N^*(v)-N^{\mathrm{EP}}(v)\right]
        \le
        1+\frac{8}{3\underline\lambda}
        +\frac{8}{3\underline\lambda}\xi(v)
        +\frac{2}{\exp(\underline\lambda^2/32)-1}.
        \label{eq:heterogeneous_phasewise_binary_expected}
    \end{align}
    Here, $b$ is held fixed, and the expectation is taken over the arrivals during phase $k$.
\end{lemma}

Lemma~\ref{lem:heterogeneous_pointwise_count} separates the effects of protection-level differences and arrival uncertainty. The first term in \eqref{eq:heterogeneous_pointwise_separation} depends on $\xi(v)$, which measures the discrepancy introduced by using $\eta_j^{\tau_k}$ in place of $\theta_j(\tau_k)$. When the number of remaining arrivals with rewards above $v$ is sufficiently close to its expectation, decisions contributing to the difference between the two counts can occur only within a window at the end of the phase whose length is bounded in terms of $\xi(v)+1$.

The second term counts periods in which this number differs from its expectation by more than the specified threshold. Hoeffding's inequality bounds the probability of each such event by a quantity that decreases exponentially with the number of remaining periods. Summing these bounds gives an expected contribution independent of the phase length, yielding the final term in \eqref{eq:heterogeneous_phasewise_binary_expected}. The complete proof is provided in Appendix~\ref{app:heterogeneous_pointwise_count}.

\subsubsection{Bounding the Loss from Protection-Level Differences}
\label{subsubsec:heterogeneous_backward_propagation}

In this subsection, we bound the contribution from protection-level differences by integrating over cutoffs. 
Recall that Lemma~\ref{lem:heterogeneous_pointwise_count} bounds the expected count gap at each cutoff $v$. 
Integrating over $v\in[0,r_{\max}]$ yields
\begin{align*}
\int_0^{r_{\max}}
\mathbb E\left[N^*(v)-N^{\mathrm{EP}}(v)\right]\,dv
\le
    \left(
        1+\frac{8}{3\underline\lambda}
        +\frac{2}{\exp(\underline\lambda^2/32)-1}
    \right)r_{\max} + \frac{8}{3\underline\lambda}\int_0^{r_{\max}}\xi(v)\,dv.
\end{align*}
The focus of this subsection is thus to found the integral in the right-hand-side.

To this end, we measure the value loss from excess or insufficient protection for each type. Consider an integer end-of-phase protection level $z\in\{0,\ldots,B\}$ for type $j$. Using capacity to accept type-$j$ demand earns $r_j$ per unit, whereas retaining it yields the benchmark continuation value. Taking $r_j$ as the opportunity cost per protected unit, $\Psi_{k+1}(z)-r_jz$ measures the net value of protecting $z$ units. The benchmark's phase-specific protection level $\theta_j(\tau_k)$ maximizes this expression. We define the \emph{protection-level value loss} as the reduction in this net value from using $z$ instead:
\begin{align}
    \mathcal L_{\Psi_{k+1}}^j(z)
    :=
    \max_{\tilde b\in\{0,\ldots,B\}}
    \left\{\Psi_{k+1}(\tilde b)-r_j\tilde b\right\}
    -\left\{\Psi_{k+1}(z)-r_jz\right\}.
    \label{eq:heterogeneous_protection_level_value_loss}
\end{align}
Since $\Psi_{k+1}$ is discretely concave, this loss can also be written as
\begin{align}
    \mathcal L_{\Psi_{k+1}}^j(z)
    =
    \sum_{u=1}^{z}\bigl(r_j-\nu_u^{(k)}\bigr)^+
    +\sum_{u=z+1}^{B}\bigl(\nu_u^{(k)}-r_j\bigr)^+.
    \label{eq:heterogeneous_protection_loss_marginals}
\end{align}
The first sum measures the value loss from protecting units whose continuation marginal values are below $r_j$; the second measures the value loss from leaving unprotected units whose continuation marginal values exceed $r_j$. Both sums are zero when $z=\theta_j(\tau_k)$, the benchmark's phase-specific protection level.

We now use this representation to bound the integral of $\xi(v)$.
Recall from \eqref{eq:heterogeneous_d_definition} that $\xi(v)$
is the largest protection discrepancy
across types at cutoff $v$.
Since these discrepancies are nonnegative, we can replace $\max$ by their sum and obtain
\begin{align}
    \xi(v)\le{}&
    \sum_{\substack{j\in\mathcal J_k\\r_j>v}}
    \left(
        \min\{\lceil\eta_j^{\tau_k}\rceil,b\}
        -\sum_{u=1}^{b}\mathbf 1\{\nu_u^{(k)}>v\}
    \right)^+
    +
    \sum_{\substack{j\in\mathcal J_k\\r_j\le v}}
    \left(
        \sum_{u=1}^{b}\mathbf 1\{\nu_u^{(k)}>v\}
        -\min\{\lceil\eta_j^{\tau_k}\rceil,b\}
    \right)^+.
    \label{eq:heterogeneous_discrepancy_sum}
\end{align}
We evaluate the integrals of these terms for each specific $j\in\mathcal J_k$.
Fixing any $j$, the excess-protection term appears when $v<r_j$,
and the insufficient-protection term appears when $v\ge r_j$.
Since the continuation marginal values are nonincreasing,
integrating over these respective ranges gives
\begin{align*}
    &\int_0^{r_j}
    \left(
        \min\{\lceil\eta_j^{\tau_k}\rceil,b\}
        -\sum_{u=1}^{b}\mathbf 1\{\nu_u^{(k)}>v\}
    \right)^+\,dv
    =
    \sum_{u=1}^{\min\{\lceil\eta_j^{\tau_k}\rceil,b\}}
    \bigl(r_j-\nu_u^{(k)}\bigr)^+,\\
    &\int_{r_j}^{r_{\max}}
    \left(
        \sum_{u=1}^{b}\mathbf 1\{\nu_u^{(k)}>v\}
        -\min\{\lceil\eta_j^{\tau_k}\rceil,b\}
    \right)^+\,dv
    =
    \sum_{u=\min\{\lceil\eta_j^{\tau_k}\rceil,b\}+1}^{b}
    \bigl(\nu_u^{(k)}-r_j\bigr)^+.
\end{align*}
The two sums on the right contain precisely the terms with $u\le b$ in \eqref{eq:heterogeneous_protection_loss_marginals}, evaluated at $z=\min\{B,\lceil\eta_j^{\tau_k}\rceil\}$. The omitted terms are nonnegative. Integrating \eqref{eq:heterogeneous_discrepancy_sum} and applying this comparison to each type therefore gives
\begin{align}
    \int_0^{r_{\max}}\xi(v)\,dv
    \le
    \sum_{j\in\mathcal J_k}
    \mathcal L_{\Psi_{k+1}}^j
    \left(\min\{B,\lceil\eta_j^{\tau_k}\rceil\}\right).
    \label{eq:heterogeneous_integrated_protection_loss}
\end{align}

We now bound the right-hand side of \eqref{eq:heterogeneous_integrated_protection_loss} by backward propagation. Fix $j\in\mathcal J_k$, and let $\kappa_j$ denote the phase in which type $j$ expires, so $e_j=\tau_{\kappa_j}$. At its expiration, $\eta_j^{e_j}=\theta_j(e_j)$, and hence $\mathcal L_{\Psi_{\kappa_j+1}}^j(\eta_j^{e_j})=\mathcal L_{\Psi_{\kappa_j+1}}^j(\theta_j(e_j))=0$.

Fix $\ell\in\{k+1,\ldots,\kappa_j\}$. We propagate the loss bound backward from the end of phase $\ell$ to the end of phase $\ell-1$. We obtain the effective end-of-phase protection level $\eta_j^{\tau_{\ell-1}}$ by adding the phase-$\ell$ demand adjustment $\sum_{s\in\mathcal T_\ell}(\sum_{i<j}\lambda_{i,s}+\frac12\lambda_{j,s})$ to $\eta_j^{\tau_\ell}$.

Let $L_\ell:=\tau_\ell-\tau_{\ell-1}$ denote the length of phase $\ell$.
The following lemma bounds the increase in protection-level value loss under this update.

\begin{lemma}[One-step bound for protection-level value loss]
\label{lem:heterogeneous_EP_one_step}
Suppose the active-interval assumption in \eqref{eq:arrival_probability_lower_bound} holds. For every phase $k$, type $j\in\mathcal J_k$, and $\ell\in\{k+1,\ldots,\kappa_j\}$,
\begin{align}
    \mathcal L_{\Psi_\ell}^j
    \left(\min\{B,\lceil\eta_j^{\tau_{\ell-1}}\rceil\}\right)
    \le
    \mathcal L_{\Psi_{\ell+1}}^j
    \left(\min\{B,\lceil\eta_j^{\tau_\ell}\rceil\}\right)
    +4r_{\max}L_\ell
    \exp\left(-\frac{\underline\lambda^2}{2}L_\ell\right).
    \label{eq:heterogeneous_EP_one_step}
\end{align}
\end{lemma}

The formal proof for the lemma is provided in Appendix~\ref{app:heterogeneous_EP_one_step}.
The key to the proof is to establish that an additional loss term is needed only when the phase-$\ell$ arrival counts deviate sufficiently from their expectations. Combining a bound on the loss for these realizations with a concentration bound on their probability gives the stated bound.

We now combine this Lemma with \eqref{eq:heterogeneous_integrated_protection_loss} and apply Lemma~\ref{lem:heterogeneous_pointwise_count} and to prove Theorem~\ref{thm:heterogeneous_EA_performance}.

\proof{Proof of Theorem~\ref{thm:heterogeneous_EA_performance}.}

Fix a phase $k$ and an entering capacity $b\in\{0,\ldots,B\}$. For each $j\in\mathcal J_k$, the protection-level value loss is zero at expiration. Applying Lemma~\ref{lem:heterogeneous_EP_one_step} successively backward from phase $\kappa_j$ to phase $k$ gives
\begin{align*}
    \mathcal L_{\Psi_{k+1}}^j
    \left(\min\{B,\lceil\eta_j^{\tau_k}\rceil\}\right)
    \le
    4r_{\max}\sum_{\ell=k+1}^{\kappa_j}
    L_\ell\exp\left(-\frac{\underline\lambda^2}{2}L_\ell\right).
\end{align*}
Summing over types and using \eqref{eq:heterogeneous_integrated_protection_loss}, we obtain
\begin{align*}
    \int_0^{r_{\max}}\xi(v)\,dv
    \le
    4r_{\max}\sum_{j\in\mathcal J_k}
    \sum_{\ell=k+1}^{\kappa_j}
    L_\ell\exp\left(-\frac{\underline\lambda^2}{2}L_\ell\right).
\end{align*}
Integrating the expected count bound in Lemma~\ref{lem:heterogeneous_pointwise_count} over $v\in[0,r_{\max}]$ and substituting this estimate yields
\begin{align}
    \int_0^{r_{\max}}
    \mathbb E\left[N^*(v)-N^{\mathrm{EP}}(v)\right]\,dv
    \le{}&
    \left(
        1+\frac{8}{3\underline\lambda}
        +\frac{2}{\exp(\underline\lambda^2/32)-1}
    \right)r_{\max}\notag\\
    &{}+\frac{32r_{\max}}{3\underline\lambda}
    \sum_{j\in\mathcal J_k}\sum_{\ell=k+1}^{\kappa_j}
    L_\ell\exp\left(-\frac{\underline\lambda^2}{2}L_\ell\right).
    \label{eq:heterogeneous_integrated_deficit}
\end{align}

We now apply this bound conditionally on the policy's entering capacity, $B_{\tau_{k-1}+1}^{\mathrm{EP}}=b$. Since the phase-$k$ arrivals are independent of this capacity, \eqref{eq:phase_aware_recursion} and \eqref{eq:heterogeneous_item_representation} give
\begin{align*}
    \Psi_k(b)
    -\mathbb E\!\left[
        R_k^{\mathrm{EP}}
        +\Psi_{k+1}(B_{\tau_k+1}^{\mathrm{EP}})
        \,\middle|\,
        B_{\tau_{k-1}+1}^{\mathrm{EP}}=b
    \right] =
    \int_0^{r_{\max}}
    \mathbb E\!\left[
        N^*(v)-N^{\mathrm{EP}}(v)
        \,\middle|\,
        B_{\tau_{k-1}+1}^{\mathrm{EP}}=b
    \right]\,dv,
\end{align*}
where the counts correspond to phase $k$ with entering capacity $b$. Note that the bound in \eqref{eq:heterogeneous_integrated_deficit} applies uniformly over $b$.

Taking expectations over the entering capacity and summing over phases cancels the intermediate continuation values. Using $\Psi_1(B)=\mathcal V_{\mathrm{PA}}$ and $\Psi_{K+1}\equiv0$, we obtain
\begin{align}
    \mathcal V_{\mathrm{PA}}
    -\mathbb E[\mathrm{Reward}(\mathrm{EP})]&=
    \sum_{k=1}^K
    \mathbb E\!\left[
        \Psi_k(B_{\tau_{k-1}+1}^{\mathrm{EP}})
        -R_k^{\mathrm{EP}}
        -\Psi_{k+1}(B_{\tau_k+1}^{\mathrm{EP}})
    \right] \nonumber \\
    & \le
    K\left(
        1+\frac{8}{3\underline\lambda}
        +\frac{2}{e^{\underline\lambda^2/32}-1}
    \right)r_{\max} +
    \frac{32r_{\max}}{3\underline\lambda}
    \sum_{k=1}^K
    \sum_{j\in\mathcal J_k}
    \sum_{\ell=k+1}^{\kappa_j}
    L_\ell e^{-\underline\lambda^2L_\ell/2}.\label{eq:vpa_bd_a}
\end{align}

For the second term, we have
\begin{align*}
    \sum_{k=1}^K\sum_{j\in\mathcal J_k}
    \sum_{\ell=k+1}^{\kappa_j}
    L_\ell e^{-\underline\lambda^2L_\ell/2}
    &\le
    \frac{2}{\underline\lambda^2}
    \sum_{k=1}^K\sum_{j\in\mathcal J_k}(\kappa_j-k) \le
    \frac{JK(K-1)}{\underline\lambda^2}.
\end{align*}
The first inequality uses $x e^{-\underline\lambda^2x/2}\le2/\underline\lambda^2$ for $x\ge0$. The second uses $\kappa_j-k\le K-k$ and the fact that at most $J$ types can arrive in each phase. Substituting this bound into \eqref{eq:vpa_bd_a} gives 
\[
    \mathcal V_{\mathrm{PA}}
    -\mathbb E[\mathrm{Reward}(\mathrm{EP})] \le K\left(1 + \frac{8}{3 \underline \lambda} + \frac{2}{\exp(\underline \lambda^2/32)-1}\right)r_{\max} +
        \frac{
            32r_{\max}JK(K-1)
        }{3\underline\lambda^3}
\]
\Halmos\endproof

\section{Numerical Experiments}
\label{sec:numerical_block_noise}

In this section, we evaluate the phase and expiration protection policies under different patterns of type disappearance and time-varying arrival probabilities. First, with known arrival probabilities, we compare their rewards with those of a budget-ratio policy (BR). This policy applies the acceptance comparison in \eqref{eq:general_acceptance} with $T$ as the protection reference time for every type and therefore a zero protection level. 
Second, we evaluate the protection policies when arrival probabilities and protection levels are estimated from sampled arrival sequences, and examine how sample size affects mean rewards and variation across independently generated data sets.

\subsection{Performance with Known Arrival Probabilities}
\label{subsec:block_noise_policy_comparison}

In this subsection, we compare the phase and expiration protection policies under known arrival probabilities and exactly computed protection levels, using the budget-ratio policy as a reference.

For each horizon $T \in \{800, 1600, 3200, 6400\}$, we consider
four types with rewards $r = (20, 12, 5, 1)$, all with activation
time $s_j = 1$ but with different expiration times.
Table~\ref{tab:block_noise_expirations} gives the two expiration
patterns, each specifying $(e_1, e_2, e_3, e_4)$, where $e_j$ is
the last period in which type~$j$ may arrive.
In Pattern~1, lower-reward types disappear earlier. Pattern~2 exchanges the expiration times of types~2 and~3.
The table also reports $\bar\lambda_j$, the mean arrival probability for type~$j$ in each period of its active interval.
\begin{table}[htbp]
    \TABLE
    {Expiration times and mean arrival probabilities\label{tab:block_noise_expirations}}
    {\begin{tabular}{lcc}
        \toprule
        Pattern & $(e_1,e_2,e_3,e_4)$ & $(\bar\lambda_1,\bar\lambda_2,\bar\lambda_3,\bar\lambda_4)$ \\
        \midrule
        1 & $(T,3T/4,T/2,T/4)$ & $(0.09,0.12,0.18,0.36)$ \\
        2 & $(T,T/2,3T/4,T/4)$ & $(0.09,0.18,0.12,0.36)$ \\
        \bottomrule
    \end{tabular}}
    {}
\end{table}

We test constant arrival probabilities following each pattern, where $\lambda_{j,t}=\bar\lambda_j$ throughout type~$j$'s active interval.
To test how the policies perform when arrival rates stay above or below their means for consecutive periods, we also let the probabilities vary across blocks within each type's active interval. Let $\rho$ denote the block length as a fraction of a type's active interval, so that each block for type~$j$ contains $\rho e_j$ periods.
We consider $\rho\in\{0.05,0.10,0.25,0.50\}$ and set the relative fluctuation amplitude to $\delta=0.25$. Independently for each type $j$ and block $h$, we draw $Z_{j,h}\in\{-1,+1\}$ with equal probabilities and set
\begin{equation}
    \lambda_{j,t}=\bar\lambda_j\left(1+\delta Z_{j,\,1+\lfloor(t-1)/(\rho e_j)\rfloor}\right),
    \qquad 1\le t\le e_j.
    \label{eq:block_noise_probabilities}
\end{equation}
The draw $Z_{j,h}$ determines whether type~$j$'s arrival probability in block $h$ is above or below $\bar\lambda_j$. Larger values of $\rho$ give longer-lasting fluctuations with the same amplitude. Conditional on the resulting probabilities, arrivals are independent across periods. Since $\mathbb E[Z_{j,h}]=0$, we have $\mathbb E[\lambda_{j,t}]=\bar\lambda_j$ while type~$j$ is active. In both cases, we set $\lambda_{j,t}=0$ for $t>e_j$. We set $B=0.225T$, the expected arrivals of types~1--2 plus half those of type~3.

For each of the two expiration patterns and four horizons, we use constant arrival probabilities ($\rho = 0$) and fluctuating probabilities at four positive values of $\rho$. This gives $2\times4\times5=40$ parameter settings in total. For each setting, we consider $S=500$ arrival-probability realizations, each specifying the arrival probabilities for every type and period. For fluctuating probabilities, we draw the full collection $\{Z_{j,h}\}$ independently for each realization. With constant probabilities, all $S$ realizations use the same arrival probabilities.

For each realization, we compute the phase-aware continuation values to obtain the phase-aware benchmark and both policies' protection levels from \eqref{eq:theta_kj_def}. 
Moreover, we generate $N=2{,}000$ independent simulated arrival sequences over the full horizon and evaluate the three online policies on the same sequences. We generate a new set of sequences for each realization, giving $SN=1{,}000{,}000$ simulated arrival sequences per parameter setting.

Let $R_{i,a,\ell}(\pi)$ denote the total reward of policy~$\pi$ on simulated arrival sequence~$\ell$ under arrival-probability realization~$a$ in parameter setting~$i$. Within each setting, we first average over the simulated arrival sequences and then over the arrival-probability realizations:
\begin{equation}
    \bar R_{i,a}(\pi)=\frac{1}{N}\sum_{\ell=1}^{N}R_{i,a,\ell}(\pi),
    \qquad \text{and} \qquad
    \widehat R_i(\pi)=\frac{1}{S}\sum_{a=1}^{S}\bar R_{i,a}(\pi).
    \label{eq:block_noise_policy_reporting}
\end{equation}
For each arrival-probability realization, we compute the phase-aware benchmark minus $\bar R_{i,a}(\pi)$ for each of the three policies, and the exact gap between the offline and phase-aware benchmarks.

Table~\ref{tab:known_probability_benchmark_gaps} reports the worst and average gaps between the phase-aware benchmark and each policy, and between the offline and phase-aware benchmarks. 
Each of the eight rows corresponds to one expiration pattern and horizon, with $B/T=0.225$ throughout. For each comparison, we report the maximum and mean gaps in the columns labeled Worst and Average, respectively. These statistics are computed over all $5S=2{,}500$ arrival-probability realizations pooled across the five values of $\rho$, and are also shown in the figures below.

\begin{table}[htbp]
    \setlength{\tabcolsep}{2pt}
    \TABLE
    {Worst and average reward gaps\label{tab:known_probability_benchmark_gaps}}
    {\begin{tabular*}{\textwidth}{@{\extracolsep{\fill}}cc*{4}{rr}@{}}
        \toprule
        Pattern & $T$ & \multicolumn{2}{c}{$\text{PA}-\text{PP}$} & \multicolumn{2}{c}{$\text{PA}-\text{EP}$} & \multicolumn{2}{c}{$\text{PA}-\text{BR}$} & \multicolumn{2}{c}{$\text{Offline}-\text{PA}$} \\
        \cmidrule(lr){3-4}\cmidrule(lr){5-6}\cmidrule(lr){7-8}\cmidrule(lr){9-10}
         & & Worst & Average & Worst & Average & Worst & Average & Worst & Average \\
        \midrule
        1 & 800 & 10.429 & 1.811 & 10.433 & 1.811 & 30.128 & 8.066 & 46.854 & 41.870 \\
         & 1600 & 15.746 & 1.694 & 15.751 & 1.694 & 43.144 & 10.576 & 66.133 & 59.676 \\
         & 3200 & 23.960 & 1.407 & 23.960 & 1.407 & 66.402 & 13.986 & 93.467 & 84.380 \\
         & 6400 & 32.211 & 1.654 & 32.211 & 1.654 & 84.349 & 19.271 & 132.047 & 119.428 \\
        \midrule
        2 & 800 & 11.718 & 2.024 & 11.715 & 2.024 & 22.937 & 8.772 & 35.558 & 26.037 \\
         & 1600 & 15.998 & 1.873 & 16.011 & 1.874 & 31.463 & 11.322 & 50.061 & 36.560 \\
         & 3200 & 21.514 & 1.582 & 21.514 & 1.582 & 43.478 & 14.970 & 70.750 & 51.518 \\
         & 6400 & 35.351 & 1.552 & 35.351 & 1.552 & 59.856 & 20.268 & 100.012 & 72.792 \\
        \bottomrule
    \end{tabular*}}
    {\NOTE{Notes.} PA denotes the phase-aware benchmark; PP, EP, and BR denote the phase protection, expiration protection, and budget-ratio policies. The Worst and Average columns report the maximum and mean total reward gaps, respectively, across the 2,500 arrival-probability realizations for each pattern and horizon. We use $B/T=0.225$ throughout.}
\end{table}

Across the eight rows of Table~\ref{tab:known_probability_benchmark_gaps}, PP and EP earn nearly identical average rewards, with average gaps to the phase-aware benchmark ranging from $1.41$ to $2.02$. EP uses four protection levels compared with ten for PP. Both protection policies have smaller average and worst observed gaps than the budget-ratio policy in every row. 
To assess the statistical significance of the protection policies' reward advantages over BR, we also compute separate $95$-percent confidence intervals for the mean differences PP--BR and EP--BR in each of the $8\times5=40$ settings, using $S=500$ paired differences per interval. All differences are statistically significant, as all intervals lie above zero.

\begin{figure}[htbp]
    \FIGURE
    {\includegraphics[width=\textwidth]{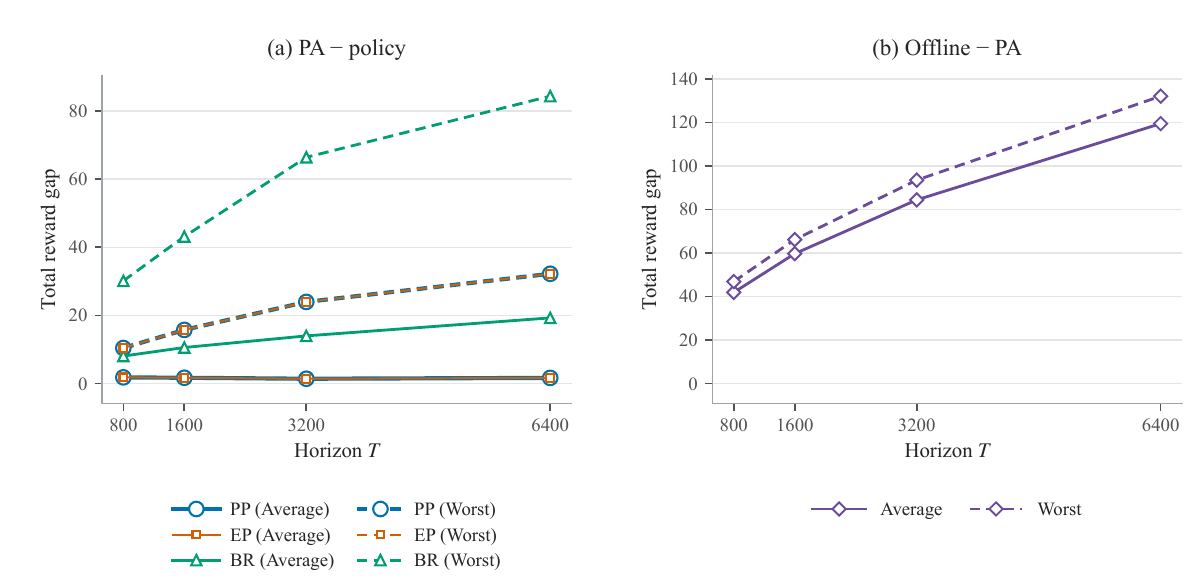}}
    {Reward gaps under Pattern~1\label{fig:block_noise_policy_case1}}
    {PA denotes the phase-aware benchmark; PP, EP, and BR denote the phase protection, expiration protection, and budget-ratio policies. All gaps are in total reward units, with $B/T=0.225$. The maximum and mean gaps are computed as in Table~\ref{tab:known_probability_benchmark_gaps}.}
\end{figure}

\begin{figure}[htbp]
    \FIGURE
    {\includegraphics[width=\textwidth]{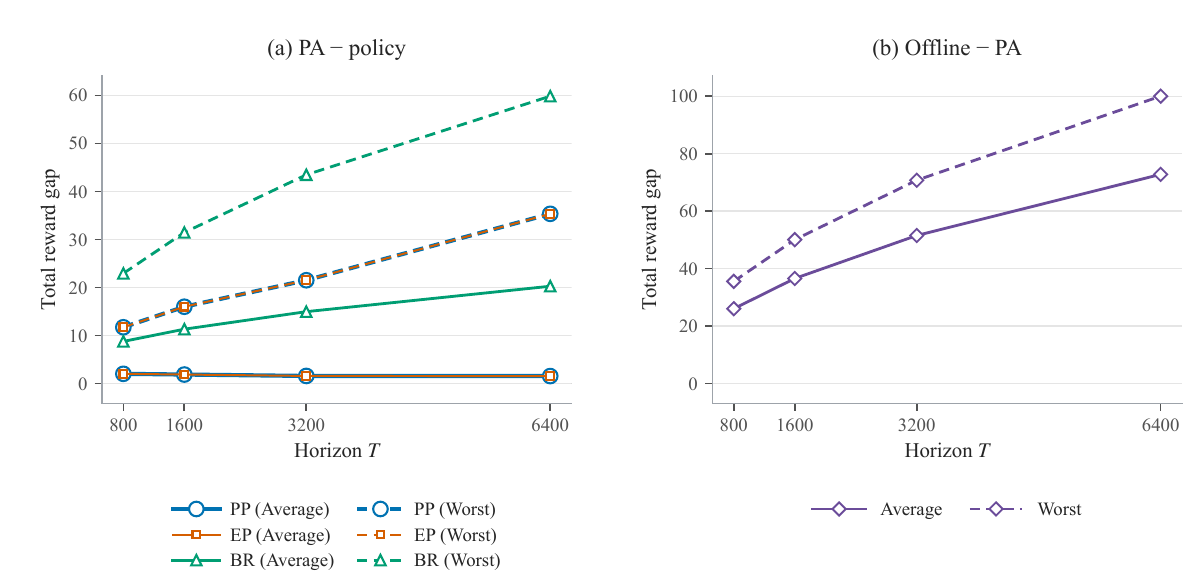}}
    {Reward gaps under Pattern~2\label{fig:block_noise_policy_case2}}
    {}
\end{figure}

Figures~\ref{fig:block_noise_policy_case1} and~\ref{fig:block_noise_policy_case2} illustrate how the worst and average gaps vary with $T$ under Patterns~1 and~2, respectively. In each figure, panel~(a) compares the phase-aware benchmark with the three policies, and panel~(b) compares the offline and phase-aware benchmarks. Solid lines show Average and dashed lines show Worst, with $B/T=0.225$ throughout.

Under Pattern~1, Figure~\ref{fig:block_noise_policy_case1}(a) shows that the protection policies' average gaps remain small as $T$ increases. 
Moreover, their advantage over the budget-ratio policy is growing with $T$. Specifically, the EP policy's mean total reward advantage over the budget-ratio policy grows from $6.26$ at $T=800$ to $17.62$ at $T=6400$.

The protection policies' average gaps also remain small under Pattern~2, where the budget-ratio policy's average gap grows from $8.77$ to $20.27$ (Figure~\ref{fig:block_noise_policy_case2}(a)). Under both patterns, the worst and average gaps between the offline and phase-aware benchmarks increase with $T$ (panel~(b) in Figures~\ref{fig:block_noise_policy_case1} and~\ref{fig:block_noise_policy_case2}).

\subsection{Performance with Estimated Arrival Probabilities}
\label{subsec:block_noise_m_study}

In this subsection, we evaluate the phase and expiration protection policies when arrival probabilities (and hence the protection levels) are estimated from sampled arrival sequences. 
We want to understand whether the EP policy is less sensitive to estimation noise, given that it uses fewer protection levels than the PP policy. To this end, we compare the policies' mean rewards, the standard deviations of their mean rewards across data sets, and the paired reward differences at different sample sizes.

We use the ten parameter settings from Subsection~\ref{subsec:block_noise_policy_comparison} with $T=1600$ and $B/T=0.225$. We follow the same procedure to generate $S=500$ arrival-probability realizations per setting and $N=2000$ evaluation sequences per realization. 
For policy construction, the rewards, expiration times, and initial capacity are known. We consider sample sizes $M\in\{2,3,5,10,20,30\}$. 
For each arrival-probability realization, we generate $Q=100$ independent data sets of arrival sequences under the same arrival probabilities to estimate the average and variability of each policy's mean reward across data sets. 
Each data set contains $30$ full-horizon arrival sequences to accommodate the largest sample size. Both policies use the first $M$ sequences of each data set, so the smaller samples are nested within the larger ones. The sequences used for estimation are independent of those used for evaluation.

We estimate continuation values by replacing the expectations in the backward recursion, including the final phase, with averages over the $M$ observed phase-arrival count vectors. Let $d_j^{(m)}(t)$ denote the indicator that a type-$j$ arrival occurs in period~$t$ of sampled arrival sequence~$m$. The empirical arrival probability is
\begin{equation}
    \widehat\lambda_{j,t}
    =\frac{1}{M}\sum_{m=1}^{M}d_j^{(m)}(t).
    \label{eq:block_noise_sample_probability}
\end{equation}
Recall that the demand adjustment in Algorithm~\ref{Algo:GeneralPolicy} is strictly positive for an arriving type because it includes half that type's arrival probability in the current period. The estimated demand sum over $[t,u]$ can be zero even when type~$j$ is active, thus we use the demand adjustment
$\max\left\{1,\,
        \sum_{\tau=t}^{u}\left(\sum_{i<j}\widehat\lambda_{i,\tau}
        +\frac12\widehat \lambda_{j,\tau}\right)
    \right\}.$
With integer capacity and protection levels, this floor preserves the requirement that at least one unit of capacity is available above the estimated protection level before accepting an arrival. We choose to add this floor to be consistent within the framework of Algorithm~\ref{Algo:GeneralPolicy}, dropping it does not meaningfully affect our numerical findings. 

For each arrival-probability realization, we use the same evaluation sequences for both policies, all data sets, and all values of $M$. Index the $10$ parameter settings by $i$. Let $R_{i,a,M,q,\ell}(\pi)$ denote the total reward on simulated arrival sequence~$\ell$ under arrival-probability realization~$a$, using the first $M$ sequences of the $q$-th data set. For the policy constructed from this data set, the mean reward over the $N$ simulated arrival sequences is
\begin{equation}
    \bar R_{i,a,M,q}(\pi)
    =\frac{1}{N}\sum_{\ell=1}^{N}R_{i,a,M,q,\ell}(\pi).
    \label{eq:block_noise_m_path_mean}
\end{equation}
For each arrival-probability realization and $M$, we repeat this construction and evaluation for the $Q=100$ independent data sets described above. This gives one mean reward per data set for each policy. We compute the mean and sample standard deviation of these $Q$ mean rewards:
\begin{align}
    \widetilde R_{i,a,M}(\pi)
    &=\frac{1}{Q}\sum_{q=1}^{Q}\bar R_{i,a,M,q}(\pi),
    \label{eq:block_noise_m_mean}\\
    s_{i,a,M}(\pi)
    &=\left[\frac{1}{Q-1}\sum_{q=1}^{Q}
    \left(\bar R_{i,a,M,q}(\pi)-\widetilde R_{i,a,M}(\pi)\right)^2\right]^{1/2}.
    \label{eq:block_noise_m_sd}
\end{align}
The standard deviation $s_{i,a,M}(\pi)$ measures how the policy's mean reward varies across data sets when the arrival probabilities and evaluation sequences are held fixed.

Table~\ref{tab:training_sample_comparison} reports results for $M=2,5,30$. Each row fixes the expiration pattern and $M$ and pools the $500$ realizations from each of the five arrival-probability settings. For each policy, we report the mean of $\widetilde R_{i,a,M}(\pi)$ over these $2{,}500$ realizations. In parentheses, we report the average construction SD $s_{i,a,M}(\pi)$ over the same realizations. The $\text{EP}-\text{PP}$ columns give the mean, minimum, and maximum of the paired differences $\widetilde R_{i,a,M}(\mathrm{EP})-\widetilde R_{i,a,M}(\mathrm{PP})$. The Positive (\%) column reports the percentage of realizations for which this difference is strictly positive.

\begin{table}[htbp]
    \setlength{\tabcolsep}{3pt}
    \TABLE
    {Sample size and policy performance\label{tab:training_sample_comparison}}
    {\begin{tabular*}{\textwidth}{@{\extracolsep{\fill}}ccrrrrrr@{}}
        \toprule
        Pattern & $M$ & \multicolumn{2}{c}{Mean reward (SD)} & \multicolumn{4}{c}{$\text{EP}-\text{PP}$} \\
        \cmidrule(lr){3-4}\cmidrule(lr){5-8}
         & & PP & EP & Mean & Minimum & Maximum & Positive (\%) \\
        \midrule
        1 & 2 & 4879.52 (19.29) & 4881.71 (20.14) & 2.190 & $-0.891$ & 10.524 & 96.96 \\
         & 5 & 4890.43 (9.29) & 4890.87 (9.83) & 0.445 & $-0.592$ & 2.312 & 79.88 \\
         & 30 & 4898.39 (1.66) & 4898.50 (1.67) & 0.111 & $-0.145$ & 0.525 & 88.32 \\
        \midrule
        2 & 2 & 4909.11 (14.98) & 4913.11 (16.81) & 3.995 & $-0.488$ & 13.316 & 99.28 \\
         & 5 & 4918.08 (8.49) & 4919.31 (9.35) & 1.223 & $-0.247$ & 3.584 & 98.04 \\
         & 30 & 4925.05 (1.66) & 4925.25 (1.66) & 0.198 & $-0.110$ & 0.582 & 99.12 \\
        \bottomrule
    \end{tabular*}}
    {\NOTE{Notes.} PP and EP denote the phase and expiration protection policies. $M$ is the number of sampled arrival sequences used to construct each policy. Each row pools $2{,}500$ arrival-probability realizations. Parentheses report the average construction SD. Positive reports the percentage of realizations in which EP's mean reward strictly exceeds PP's. Reward statistics are in total reward units, with $T=1600$ and $B/T=0.225$.}
\end{table}

\begin{figure}[htbp]
    \FIGURE
    {\includegraphics[width=\textwidth]{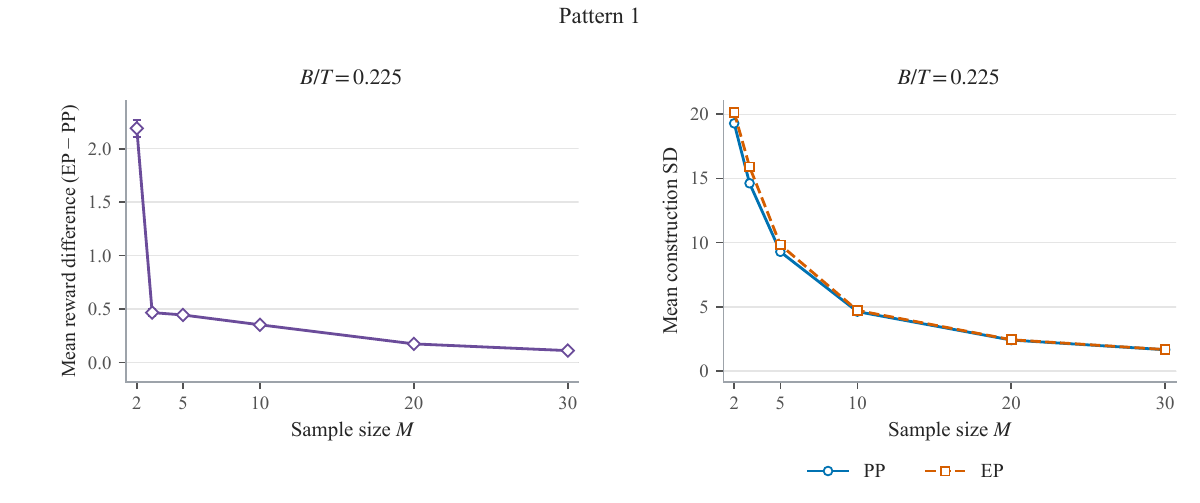}}
    {Effect of sample size under Pattern~1\label{fig:block_noise_m_case1}}
    {PP and EP denote the phase and expiration protection policies. The left panel shows the mean reward difference EP$-$PP; the right panel shows Mean construction SD. Averaging and units follow Table~\ref{tab:training_sample_comparison}. Error bars extend $1.96$ standard errors in each direction. Within each setting, the standard error is the sample SD of the 500 paired reward differences $\widetilde R_{i,a,M}(\mathrm{EP})-\widetilde R_{i,a,M}(\mathrm{PP})$, divided by $\sqrt{500}$. The standard error of the average across the five independent settings is the square root of the sum of their squared standard errors, divided by five.}
\end{figure}

\begin{figure}[htbp]
    \FIGURE
    {\includegraphics[width=\textwidth]{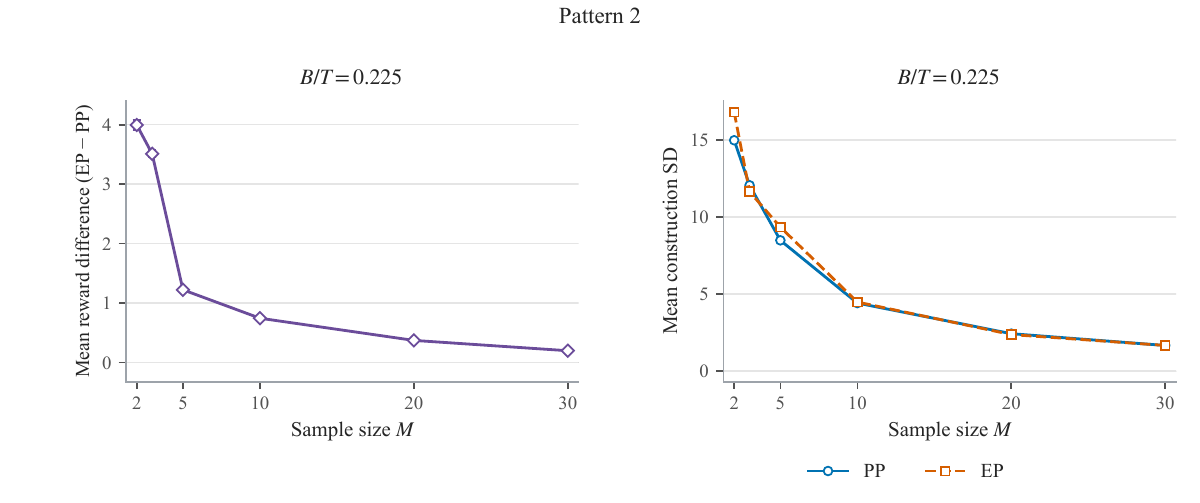}}
    {Effect of sample size under Pattern~2\label{fig:block_noise_m_case2}}
    {Mean reward differences and mean construction SDs against sample size $M$. Averaging, confidence intervals for reward differences, units, and panel layout follow Figure~\ref{fig:block_noise_m_case1}.}
\end{figure}

The mean EP--PP difference in Table~\ref{tab:training_sample_comparison}
is positive in all six rows. From $M=2$ to $M=30$, the mean advantage
falls from $2.190$ to $0.111$ under Pattern~1 and from $3.995$ to $0.198$
under Pattern~2. At the sample sizes shown, EP earns more in $79.88\%$
to $96.96\%$ of realizations under Pattern~1 and in $98.04\%$ to $99.28\%$ under Pattern~2. Under both patterns, each policy's mean reward rises and its mean construction SD falls as $M$ increases.

The simulation results suggest that the EP policy's use of fewer protection levels may indeed explain its superior performance under noisy arrival-probability estimation. Specifically, EP's mean reward advantage is larger at $M=2$ than at higher values of $M$, because fewer samples yield noisier estimates.
Recall that the EP policy uses four protection levels (one per type), whereas the PP policy uses ten (one for each active phase--type pair). Because continuation values are estimated from sampled phase-arrival counts, sampling noise can move an estimated continuation marginal value across the corresponding reward $r_j$, producing an erroneous threshold. Relying on more thresholds computed from noisy estimates creates more points at which such errors can lead to poor decisions (akin to overfitting), hurting the PP policy more than the EP policy. Figures~\ref{fig:block_noise_m_case1} and~\ref{fig:block_noise_m_case2} extend the comparison to all six values of $M$. As $M$ increases, the left panels show that the difference in mean reward narrows, and the right panels show that both policies' mean construction SDs decrease.

\section{Conclusion}\label{sec:conclusions}

In this work, we consider online selection with demand types that are available during type-specific intervals and disappear when those intervals end. Such arrival patterns are often studied in classical protection-level control. In our setting, however, the standard hindsight offline benchmark commonly used in the constant-regret literature can be too loose for establishing a constant optimality gap. We show that the offline benchmark value can exceed the optimal
online value by $\Omega(\sqrt T)$, even in a two-type instance. Hence, we introduce a tighter phase-aware benchmark that partitions the horizon into consecutive phases ending at the distinct type-expiration times. 

We propose two budget-ratio protection policies based on different protection reference times: the phase protection policy uses the end of the current phase, whereas the EP policy uses each type's expiration time and requires only one protection level per type. Both policies achieve an $O(1)$ expected reward gap relative to the phase-aware benchmark and hence an $O(1)$ optimality gap, uniformly in the horizon and initial capacity under our model assumptions. For the PP policy, we compare its decisions with optimal benchmark decisions period by period and use concentration inequalities to bound the expected loss. 
The analysis for the EP policy is complicated as its
effective end-of-phase protection levels need not match the benchmark's phase-specific protection levels, and its acceptance decisions need not preserve reward-priority ordering. Therefore, for the EP policy, we instead compare earned rewards and the benchmark continuation value of remaining capacity over each phase, separating arrival uncertainty from protection-level value losses. We control these two contributions using concentration inequalities and backward propagation, respectively.

Our numerical experiments evaluate the two policies under both known arrival probabilities and probabilities estimated from sampled arrival sequences. In the tested settings with known arrival probabilities, PP and EP earn nearly identical mean rewards, and both earn higher mean rewards than the budget-ratio policy. When arrival probabilities are estimated from sampled arrival sequences, EP earns higher mean rewards than PP, with the advantage narrowing as sample size increases.

\bibliographystyle{informs2014}

\bibliography{Resource_allocation}

\ECSwitch

\newcommand*{\ECHyperrefAnchor}[1]{%
  \expandafter\def\csname theH#1\endcsname{EC.\arabic{#1}}%
}
\ECHyperrefAnchor{section}
\ECHyperrefAnchor{equation}
\ECHyperrefAnchor{theorem}
\ECHyperrefAnchor{lemma}
\ECHyperrefAnchor{proposition}
\ECHyperrefAnchor{corollary}
\ECHyperrefAnchor{claim}
\ECHyperrefAnchor{conjecture}
\ECHyperrefAnchor{hypothesis}
\ECHyperrefAnchor{assumption}
\ECHyperrefAnchor{remark}
\ECHyperrefAnchor{example}
\ECHyperrefAnchor{problem}
\ECHyperrefAnchor{definition}
\ECHyperrefAnchor{question}
\ECHyperrefAnchor{answer}
\ECHyperrefAnchor{exercise}
\ECHyperrefAnchor{figure}
\ECHyperrefAnchor{table}

\ECHead{Proofs and Additional Numerical Results}

\section{Proof of Proposition~\ref{prop:sqrt_T}}

\proof{Proof.}
Let $n=T/2$, and let $N_1$ and $N_2$ be the arrival counts of types~1 and~2,
respectively. The counts are independent, each with distribution
$\mathrm{Bin}(n,1/2)$ and mean $\mu=n/2$. Since $B=3n/4$, we have
$(B+\mu)/2=5T/16$.

An optimal offline policy accepts $\min\{B,N_1+N_2\}$ arrivals in total.
It gives priority to type~1, so $\min\{B,N_1\}$ of the accepted arrivals
are of type~1. Each accepted arrival earns a base reward of $1/2$, with
an additional $1/2$ for type~1. Therefore,
\begin{align*}
    \mathcal V_{\mathrm{off}}
    &= \tfrac12\mathbb E\!\left[\min\{B,N_1+N_2\}+\min\{B,N_1\}\right]\\
    &= \tfrac{5}{16}T
       -\tfrac12\mathbb E\!\left[(N_1-B)^++(B-N_1-N_2)^+\right].
\end{align*}
The first positive part counts type-1 arrivals in excess of capacity,
and the second counts unused units of capacity. For the first term, $B$ is
$n/4$ above the mean of $N_1$. For the second, $B$ is $n/4$ below
the mean of $N_1+N_2$.

For $a>0$, we have $(x-a)^+\le x^2/a$, since $x>a$ implies $x^2/a\ge x$.
Thus, for a centered random variable, the expected excess above $a$ is
at most its variance divided by $a$. Each count has variance $n/4$, and
independence gives variance $n/2$ for their sum. Applying the bound with
$a=n/4$ gives
\begin{align*}
    \mathbb E[(N_1-B)^+]
    &= \mathbb E[((N_1-\mu)-n/4)^+]
    \le \frac{\operatorname{Var}(N_1)}{n/4}=1,\\
    \mathbb E[(B-N_1-N_2)^+]
    &= \mathbb E[((n-N_1-N_2)-n/4)^+]
    \le \frac{\operatorname{Var}(N_1+N_2)}{n/4}=2.
\end{align*}
Substituting these bounds gives $\mathcal V_{\mathrm{off}}\ge 5T/16-3/2$.

For a feasible online policy $\pi$, let $b'=B_{n+1}^\pi$ be the capacity
carried into the second half. The policy earns $\tfrac12(B-b')$ from
the type-2 arrivals it accepts in the first half. In the second half,
it earns at most $\min\{b',N_1\}$. Using
$2\min\{b',N_1\}=b'+N_1-|N_1-b'|$, we obtain
\begin{align*}
    \mathbb E[\mathrm{Reward}(\pi)]
    &\le \mathbb E\!\left[\tfrac12(B-b')+\min\{b',N_1\}\right]\\
    &= \tfrac{5}{16}T-\tfrac12\mathbb E|N_1-b'|.
\end{align*}
The policy chooses $b'$ using information from the first half. Since arrivals
are independent across periods, $b'$ is independent of $N_1$.
The function $b\mapsto\mathbb E|N_1-b|$ is convex and, by
symmetry of $N_1$, symmetric about $\mu$. It is therefore minimized
at $\mu$. Conditioning on $b'$ gives
$\mathbb E|N_1-b'|\ge \mathbb E|N_1-\mu|$.

To bound this mean absolute deviation, write $X=N_1-\mu$ as the sum
of $n$ independent symmetric $\pm 1/2$ variables. Its second moment
is $n/4$. In the fourth-moment expansion, only terms with even powers
of each variable survive, giving
$\mathbb E[X^4]=(n+6\binom{n}{2})/16\le 3n^2/16$.
Applying H\"older's inequality to $|X|^{2/3}|X|^{4/3}$ gives
$\mathbb E[X^2]\le(\mathbb E|X|)^{2/3}(\mathbb E[X^4])^{1/3}$.
Rearranging and using $\mathbb E[X^2]=n/4$, we obtain
\begin{align*}
    \mathbb E|N_1-\mu|
    \ge \frac{(\mathbb E[X^2])^{3/2}}{(\mathbb E[X^4])^{1/2}}
    \ge \sqrt{n/12}.
\end{align*}
Since $n=T/2$, the online reward is at most $5T/16-\sqrt{T}/(4\sqrt{6})$.
Together with the offline bound, this proves the claim with $M=3/2$.
\Halmos
\endproof

\section{Proofs for the Phase-Aware Benchmark}

In this section, we prove the validity of the phase-aware benchmark and characterize its optimal allocation.
First, we prove Lemma~\ref{lmm:benchmark_validity}.
Subsection~\ref{app:direct_structural_properties} then establishes discrete concavity and bounded marginals of the continuation values, the threshold and nesting properties of the protection levels, and marginal optimality conditions for the benchmark's inner problem.
Finally, Subsection~\ref{app:direct_inner_optimizer} uses these properties to prove Proposition~\ref{prop:direct_inner_optimizer}.

\subsection{Proof of Lemma~\ref{lmm:benchmark_validity}}

\proof{Proof.}
For $t\in[T]$ and $b\in\{0,\ldots,B\}$, let $V_t(b)$ denote the optimal expected reward from periods $t,\ldots,T$ over feasible online policies, with $b$ units of capacity before the period-$t$ arrival outcome is observed. 
We prove by backward induction that, for every $k\in\{1,\ldots,K+1\}$ and $b\in\{0,\ldots,B\}$,
\begin{align}
    V_{\tau_{k-1}+1}(b)
    \le
    \Psi_k(b)
    \le
    \mathbb E\left[
        P\left(b,\sum_{\ell=k}^K\boldsymbol A^{(\ell)}\right)
    \right].
    \label{eq:phase_aware_induction_bounds}
\end{align}
When $k=K+1$, 
\eqref{eq:phase_aware_induction_bounds} trivially holds as $V_{T+1}\equiv 0$ and 
$\Psi_{K+1} \equiv 0$.
Fix $k\in[K]$ and suppose that \eqref{eq:phase_aware_induction_bounds} holds for phase $k+1$ and every capacity level.

We first prove $V_{\tau_{k-1}+1}(b)\le\Psi_k(b)$. Consider any feasible online policy $\pi$ starting phase $k$ with capacity $b$. Let $R_k^\pi$ denote its reward during phase $k$, and let $b'$ denote its remaining capacity at the end of that phase. The policy accepts $b-b'$ arrivals during phase $k$, so, for every realization,
\begin{align*}
    R_k^\pi
    \le
    P(b-b',\boldsymbol A^{(k)}).
\end{align*} 
Because future arrivals are independent of the history observed up to time $\tau_k$, the policy's conditional expected reward after phase $k$ is at most $V_{\tau_k+1}(b')$. By the induction hypothesis, this is at most $\Psi_{k+1}(b')$. Therefore, the policy's expected reward from phase $k$ onward is at most
\begin{align*}
    \mathbb E\left[
        R_k^\pi+V_{\tau_k+1}(b')
    \right]
    &\le
    \mathbb E\left[
        P(b-b',\boldsymbol A^{(k)})
        +\Psi_{k+1}(b')
    \right]\\
    &\le
    \mathbb E\left[
        \max_{\tilde b\in\{0,\ldots,b\}}
        \left\{
            P(b-\tilde b,\boldsymbol A^{(k)})
            +\Psi_{k+1}(\tilde b)
        \right\}
    \right]\\
    &=\Psi_k(b).
\end{align*}
Taking the supremum over feasible online policies gives $V_{\tau_{k-1}+1}(b)\le\Psi_k(b)$.

We next prove the second inequality in \eqref{eq:phase_aware_induction_bounds}. For any realization of the arrivals from phase $k$ onward,
\begin{align}
    &\max_{\tilde b\in\{0,\ldots,b\}}
    \left\{
        P(b-\tilde b,\boldsymbol A^{(k)})
        +
        P\left(\tilde b,\sum_{\ell=k+1}^K\boldsymbol A^{(\ell)}\right)
    \right\}
    =
    P\left(b,\sum_{\ell=k}^K\boldsymbol A^{(\ell)}\right).
    \label{eq:offline_phase_split}
\end{align}
Indeed, $\tilde b$ allocates capacity to the arrivals after phase $k$, while $b-\tilde b$ allocates capacity to the current phase. Any such allocation is feasible for the offline problem. Conversely, an optimal offline allocation induces a feasible division of capacity between these two groups of arrivals.

Using the induction hypothesis in the benchmark recursion \eqref{eq:phase_aware_recursion} gives
\begin{align*}
    \Psi_k(b)
    &\le
    \mathbb E_{\boldsymbol A^{(k)}}\left[
        \max_{\tilde b\in\{0,\ldots,b\}}
        \left\{
            P(b-\tilde b,\boldsymbol A^{(k)})
            +
            \mathbb E\left[
                P\left(
                    \tilde b,\sum_{\ell=k+1}^K\boldsymbol A^{(\ell)}
                \right)
            \right]
        \right\}
    \right]\\
    &\le
    \mathbb E\left[
        \max_{\tilde b\in\{0,\ldots,b\}}
        \left\{
            P(b-\tilde b,\boldsymbol A^{(k)})
            +
            P\left(
                \tilde b,\sum_{\ell=k+1}^K\boldsymbol A^{(\ell)}
            \right)
        \right\}
    \right]\\
    &=
    \mathbb E\left[
        P\left(b,\sum_{\ell=k}^K\boldsymbol A^{(\ell)}\right)
    \right].
\end{align*}
In the first line, the inner expectation is over future phase arrivals, which are independent of $\boldsymbol A^{(k)}$. The second inequality uses the fact that the maximum of expectations is no larger than the expectation of the maximum. The final equality follows from \eqref{eq:offline_phase_split}. This completes the induction.

Finally, setting $k=1$ and $b=B$ in
\eqref{eq:phase_aware_induction_bounds} yields
\begin{align*}
    V^\star
    =
    V_1(B)
    \le
    \Psi_1(B)
    =
    \mathcal V_{\mathrm{PA}}
    \le
    \mathbb E[P(B,\boldsymbol A)]
    =
    \mathcal V_{\mathrm{off}}.
\end{align*}
\Halmos
\endproof

\subsection{Structural properties of continuation values}
\label{app:direct_structural_properties}

For a function $g$ on consecutive integers, write
\begin{align*}
    g'_-(m):=g(m)-g(m-1),
    \qquad
    g'_+(m):=g(m+1)-g(m)
\end{align*}
whenever both arguments are in the domain. For $P$, these differences are
taken with respect to capacity, with the arrival counts held fixed.

\begin{lemma}[Discrete concavity of continuation values]
\label{lem:psi_discrete_concavity}
    For every $B\geq 1$, $k\in\{1,\ldots,K+1\}$, the function $\Psi_k:\{0,\dots,B\}\to\mathbb R$ is nondecreasing and discretely concave.
    Moreover, for every $b\in\{1,\dots,B\}$, $0\le\Psi'_{k,-}(b)\le r_{\max}.$
\end{lemma}
\proof{Proof.}
We prove by induction on $k$. Clearly, the statement holds for $k =K+1$ as $\Psi_{K+1}\equiv0$. Assume it holds for $\Psi_{k+1}$ and fix a realization
of $\boldsymbol A^{(k)}$. The marginal values of capacity in the current
phase and in later phases are, respectively,
\begin{align*}
    \alpha_m:=P'_-(m,\boldsymbol A^{(k)}),
    \qquad
    \beta_m:=\Psi'_{k+1,-}(m),
    \qquad m=1,\ldots,B.
\end{align*}
The marginal $\alpha_m$ is the $m$th largest realized arrival reward,
or zero if there are fewer than $m$ arrivals. Thus the $\alpha_m$ form
a nonincreasing sequence in $[0,r_{\max}]$. The $\beta_m$ have the same
properties by the induction hypothesis.

Both value functions are zero at zero capacity. The value of any split
$\tilde b\in\{0,\ldots,b\}$ in \eqref{eq:phase_aware_recursion} is the
sum of the corresponding marginals,
\begin{align*}
    P(b-\tilde b,\boldsymbol A^{(k)})+\Psi_{k+1}(\tilde b)
    =
    \sum_{m=1}^{b-\tilde b}\alpha_m+\sum_{m=1}^{\tilde b}\beta_m .
\end{align*}
Merge the two lists into a nonincreasing sequence
$\gamma_1,\ldots,\gamma_{2B}$, preserving the order of tied entries
within each list. Every split selects $b$ entries in total, so its value
is at most $\sum_{m=1}^b\gamma_m$. The first $b$ entries of the merged
list attain this bound. Since the merge preserves the order within each
original list, these entries take an initial segment from each. If
$\tilde b$ of them come from the continuation list, the split allocates
$b-\tilde b$ units to the current phase and carries $\tilde b$ forward.

The inner maximum in \eqref{eq:phase_aware_recursion} is therefore
$\sum_{m=1}^b\gamma_m$, with the empty sum equal to zero at $b=0$.
For $b=1,\ldots,B$, the difference in the maxima at capacities $b$ and
$b-1$ is $\gamma_b$. Taking expectations gives
\begin{align*}
    \Psi_k(b)-\Psi_k(b-1)
    =\mathbb E_{\boldsymbol A^{(k)}}[\gamma_b].
\end{align*}
For each realization, $\gamma_b$ lies in $[0,r_{\max}]$ and is
nonincreasing in $b$. The expected marginals inherit these properties,
which completes the induction.
\Halmos
\endproof

\begin{corollary}[Characterization and nesting of protection levels]
\label{cor:direct_protection_levels}
    For every phase $k\in[K]$ and type $j\in\mathcal J_k$, the protection level $\theta_j(\tau_k)$ is well defined and satisfies
    \begin{align*}
        \Psi'_{k+1,-}(m)
        &\ge r_j,
        \quad m\in\{1,\ldots,\theta_j(\tau_k)\},\\
        \Psi'_{k+1,-}(m)
        &< r_j,
        \quad m\in\{\theta_j(\tau_k)+1,\ldots,B\}.
    \end{align*}
Moreover, for any $i,j\in\mathcal J_k$ with $i<j$,
    \begin{align*}
        \theta_i(\tau_k)\le\theta_j(\tau_k).
    \end{align*}
\end{corollary}
\proof{Proof.}
Fix $k$ and $j$. The set in \eqref{eq:theta_kj_def} is finite and
contains $0$, so $\theta_j(\tau_k)$ is well defined.
For $\theta_j(\tau_k)<m\le B$, the marginal is strictly below $r_j$,
since a marginal of at least $r_j$ at $m$ would make the maximum in
\eqref{eq:theta_kj_def} at least $m$. When $\theta_j(\tau_k)=0$, this
covers every $m\in\{1,\ldots,B\}$ and the first range is empty.

If $\theta_j(\tau_k)>0$, the marginal at that index is at least $r_j$
by definition. By Lemma~\ref{lem:psi_discrete_concavity}, the marginals
are nonincreasing in $m$, so in this case
\begin{align*}
    \Psi'_{k+1,-}(m)
    \ge \Psi'_{k+1,-}\bigl(\theta_j(\tau_k)\bigr)
    \ge r_j,
    \qquad 1\le m\le\theta_j(\tau_k).
\end{align*}

If $i<j$, then
$r_i\ge r_j$, so every index that meets the threshold for type $i$ also
meets the threshold for type $j$. The maxima in \eqref{eq:theta_kj_def}
then satisfy $\theta_i(\tau_k)\le\theta_j(\tau_k)$.
\Halmos
\endproof

\begin{lemma}[Marginal optimality for the benchmark inner problem]
\label{lem:multi_marginal_optimality}
    Fix a phase $k$, an inherited capacity $b\in\{0,\dots,B\}$, and a realized arrival vector $\boldsymbol A^{(k)}$. Define
    \begin{align*}
        f_k(\tilde b,\boldsymbol A^{(k)};b)
        :=
        P(b-\tilde b,\boldsymbol A^{(k)})
        +
        \Psi_{k+1}(\tilde b),
        \qquad
        \tilde b\in\{0,1,\dots,b\}.
    \end{align*}
    Then $\tilde b\mapsto f_k(\tilde b,\boldsymbol A^{(k)};b)$ is discretely concave.

    For $b\ge2$, an interior carry-over level
    $\tilde b\in\{1,\dots,b-1\}$ is optimal if and only if
    \begin{align}
        -P'_+(b-\tilde b,\boldsymbol A^{(k)})
        +
        \Psi'_{k+1,-}(\tilde b)
        \ge
        0
        \ge
        -P'_-(b-\tilde b,\boldsymbol A^{(k)})
        +
        \Psi'_{k+1,+}(\tilde b).
        \label{eq:multi_inner_optimality}
    \end{align}
    For $b\ge1$, the lower boundary $\tilde b=0$ is optimal if and only if $-P'_-(b,\boldsymbol A^{(k)}) + \Psi'_{k+1,+}(0) \le 0,$ and the upper boundary $\tilde b=b$ is optimal if and only if $-P'_+(0,\boldsymbol A^{(k)}) +  \Psi'_{k+1,-}(b) \ge 0.$
    When $b=0$, $\tilde b=0$ is the unique feasible carry-over capacity.
\end{lemma}
\proof{Proof.}
When $b=0$, the feasible set is the singleton $\{0\}$. For $b\ge1$,
write $f(m)=f_k(m,\boldsymbol A^{(k)};b)$. Each unit carried forward is
taken from the current phase, so the marginals of $f$ are
\begin{align*}
    f'_-(m)
    &= -P'_+(b-m,\boldsymbol A^{(k)})+\Psi'_{k+1,-}(m),
    &&1\le m\le b,\\
    f'_+(m)
    &= -P'_-(b-m,\boldsymbol A^{(k)})+\Psi'_{k+1,+}(m),
    &&0\le m<b.
\end{align*}
As $m$ increases, the current-phase capacity $b-m$ falls. Because the
current-phase marginals are nonincreasing in capacity, the reward given
up, $P'_+(b-m,\boldsymbol A^{(k)})$, is nondecreasing in $m$. The
continuation marginal $\Psi'_{k+1,-}(m)$ is nonincreasing by
Lemma~\ref{lem:psi_discrete_concavity}. Thus $f'_-(m)$ is nonincreasing,
proving discrete concavity.

Comparison with the feasible neighbors shows that an optimal $\tilde b$
must satisfy $f'_-(\tilde b)\ge0$ if $\tilde b>0$ and
$f'_+(\tilde b)\le0$ if $\tilde b<b$. These conditions also suffice.
Since $f'_-(m)$ is nonincreasing, a nonnegative $f'_-(\tilde b)$ makes
all earlier increments nonnegative. Likewise, a nonpositive
$f'_+(\tilde b)=f'_-(\tilde b+1)$ makes all later increments nonpositive.
Thus $f$ is nondecreasing up to $\tilde b$ and nonincreasing thereafter,
with one of those ranges empty at an endpoint. Substituting the marginal
formulas gives \eqref{eq:multi_inner_optimality} in the interior and the
stated one-sided conditions at the boundaries.
\Halmos
\endproof

\subsection{Proof of Proposition~\ref{prop:direct_inner_optimizer}}
\label{app:direct_inner_optimizer}

\proof{Proof.}
Fix a phase $k$, an inherited capacity $b$, and a realized arrival vector
$\boldsymbol A^{(k)}$. List the realized arrivals in increasing type-index
order, so that their rewards are nonincreasing. With $q$ units allocated
to the current phase, the benchmark's inner objective is
\begin{align*}
    F(q):=P(q,\boldsymbol A^{(k)})+\Psi_{k+1}(b-q),
    \qquad q=0,\ldots,b.
\end{align*}
This is the objective in Lemma~\ref{lem:multi_marginal_optimality}
indexed by $q=b-\tilde b$, so $F$ is discretely concave. Let $q^*$ be the
number of strictly positive increments of $F$. Because the increments
are nonincreasing, $F$ increases up to $q^*$ and is nonincreasing
thereafter. Thus $q^*$ maximizes $F$, with ties resolved in favor of
carrying capacity forward.

For $1\le m\le b$, if the $m$th arrival in the list has type $j$, giving
the current phase its $m$th unit adds reward $r_j$ and reduces the
carry-over capacity from $b-m+1$ to $b-m$. Hence
\begin{align*}
    F(m)-F(m-1)=r_j-\Psi'_{k+1,-}(b-m+1).
\end{align*}
If the list contains fewer than $m$ arrivals, the current-phase reward
is already constant and $F(m)-F(m-1)\le 0$ by monotonicity of $\Psi_{k+1}$.
Thus $q^*$ is at most the number of realized arrivals.

By Corollary~\ref{cor:direct_protection_levels}, for a type-$j$ arrival
at rank $m\le b$,
\begin{align*}
    F(m)-F(m-1)>0
    \quad\Longleftrightarrow\quad b-m+1>\theta_j(\tau_k)
    \quad\Longleftrightarrow\quad m\le b-\theta_j(\tau_k).
\end{align*}
Since $\theta_j(\tau_k)\ge0$, ranks above $b$ also fail the last inequality.
Thus the first $q^*$ arrivals are exactly those whose ranks satisfy the
bound for their type. Let
$S_{<j}:=\sum_{i\in\mathcal J_k,\,i<j} A_i^{(k)}$ count the realized arrivals
of types preceding $j$. The type-$j$ arrivals occupy ranks
$S_{<j}+1,\ldots,S_{<j}+A_j^{(k)}$. The number of these ranks at most
$b-\theta_j(\tau_k)$ is
\begin{align*}
    \min\left\{A_j^{(k)},\left(b-\theta_j(\tau_k)-S_{<j}\right)^+\right\}
    = x_j^*
\end{align*}
by \eqref{eq:xstar_direct}. These counts select the $q^*$ highest-reward
arrivals and hence attain $P(q^*,\boldsymbol A^{(k)})$. The resulting
carry-over capacity $b^*=b-q^*$ is optimal because $q^*$ maximizes $F$.
\Halmos
\endproof

\section{Proof of Proposition~\ref{prop:direct_tracking_error}}
\label{app:direct_tracking_error}

\proof{Proof.}
Fix a phase $k\in[K]$, a period $t\in\mathcal T_k$, and a type $j\in\mathcal J_k$ with $\lambda_{j,t}>0$. Let $H_{k,t}:=\tau_k-t+1$ denote the number of periods remaining in phase $k$, including period $t$. Suppose that $H_{k,t}\ge4/\underline\lambda$, and condition on $d_j(t)=1$. For every $i\in\mathcal J_k$, let $\mu_i^{k,t}:=\mathbb E[\widehat A_i^{(k)}(t)]=\sum_{s=t}^{\tau_k}\lambda_{i,s}$. Here, $\mu_i^{k,t}$ is the unconditional expectation used in
the policy threshold. Since $\lambda_{j,t}>0$ and $j\in\mathcal J_k$, we have $s_j\le t\le\tau_k\le e_j$. The active-interval assumption gives $\lambda_{j,s}\ge\underline\lambda$ for $s=t,\ldots,\tau_k$, so $\mu_j^{k,t}\ge\underline\lambda H_{k,t}\ge4$.

Define the residual capacities based on expected and realized arrival counts, respectively, by
\begin{align*}
    Q_j^{k,t}
    &:=B_t^{\mathrm{PP}}-\theta_j(\tau_k)
    -\sum_{\substack{i\in\mathcal J_k\\i<j}}\mu_i^{k,t},
    &\widehat Q_j^{k,t}
    &:=B_t^{\mathrm{PP}}-\theta_j(\tau_k)
    -\sum_{\substack{i\in\mathcal J_k\\i<j}}\widehat A_i^{(k)}(t).
\end{align*}
Since $\frac12\mu_j^{k,t}\ge2$, the threshold inequality in \eqref{eq:general_acceptance} already implies $B_t^{\mathrm{PP}}\ge1$. Thus, conditional on $d_j(t)=1$, the policy accepts if and only if $Q_j^{k,t}\ge\frac12\mu_j^{k,t}$. 

For satisfying actions, we claim that
\begin{align*}
    \widehat Q_j^{k,t}\ge1
    &\quad\Longrightarrow\quad
    1\in\mathcal S_{j,t}^{(k)}(B_t^{\mathrm{PP}}),\\
    \widehat Q_j^{k,t}\le\widehat A_j^{(k)}(t)-1
    &\quad\Longrightarrow\quad
    0\in\mathcal S_{j,t}^{(k)}(B_t^{\mathrm{PP}}).
\end{align*}
To prove these implications, apply Proposition~\ref{prop:direct_inner_optimizer}
to the benchmark problem in \eqref{eq:direct_W_EX} with capacity
$B_t^{\mathrm{PP}}$ and remaining arrival counts $\widehat{\boldsymbol A}^{(k)}(t)$.
The number of type-$j$ arrivals accepted in the resulting optimal allocation is
\begin{align*}
    x_j^*
    \bigl(B_t^{\mathrm{PP}},\widehat{\boldsymbol A}^{(k)}(t)\bigr)
    =
    \min\left\{
        \widehat A_j^{(k)}(t),
        \left(\widehat Q_j^{k,t}\right)^+
    \right\}.
\end{align*}
Under $d_j(t)=1$, we have $\widehat A_j^{(k)}(t)\ge1$. If $\widehat Q_j^{k,t}\ge1$, this optimal allocation accepts at least one type-$j$ arrival. Since arrivals of the same type have the same reward and capacity requirement, we can choose the current arrival as one of them without changing the objective value. Thus, the acceptance branch in \eqref{eq:satisfying_equation} attains the maximum, and acceptance is satisfying.
If $\widehat Q_j^{k,t}\le\widehat A_j^{(k)}(t)-1$, the optimal allocation accepts at most $\widehat A_j^{(k)}(t)-1$ type-$j$ arrivals. These acceptances can all be assigned to later arrivals of that type, leaving the current arrival rejected without changing the objective value. The rejection branch in \eqref{eq:satisfying_equation} therefore attains the maximum, and rejection is satisfying. These conditions are sufficient and need not be necessary when the benchmark optimizer is nonunique.

By definition, $Q_j^{k,t}-\widehat Q_j^{k,t}=\sum_{i\in\mathcal J_k,\,i<j}(\widehat A_i^{(k)}(t)-\mu_i^{k,t})$. We consider two types of tracking errors.

\paragraph{Case 1: Wrong rejection.}
Suppose the policy rejects but rejection is not satisfying. Then $Q_j^{k,t}<\frac12\mu_j^{k,t}$ and $\widehat Q_j^{k,t}>\widehat A_j^{(k)}(t)-1$, which imply
\begin{align}
    \sum_{\substack{i\in\mathcal J_k\\ i\le j}}
    \left(
        \widehat A_i^{(k)}(t)
        -
        \mathbb E[\widehat A_i^{(k)}(t)]
    \right)
    <
    -\frac{1}{2}\mu_j^{k,t}+1
    \le
    -\frac14\mu_j^{k,t}.
    \label{eq:wrong_rejection_deviation}
\end{align}

For $s=t,\ldots,\tau_k$, let $X_s:=\sum_{i\in\mathcal J_k,\,i\le j}d_i(s)$ and $p_{\le j,s}:=\sum_{i\in\mathcal J_k,\,i\le j}\lambda_{i,s}$. Conditional on $d_j(t)=1$, we have $X_t=1$, so the current-period centered term is $1-p_{\le j,t}\ge0$. Thus, \eqref{eq:wrong_rejection_deviation} implies $\sum_{s=t+1}^{\tau_k}(X_s-p_{\le j,s})<-\frac14\mu_j^{k,t}$. Independence across periods ensures that the future $X_s$ remain independent Bernoulli random variables with means $p_{\le j,s}$. Hoeffding's inequality gives
\begin{align*}
    \mathbb P
    \left(
        \text{wrong rejection}
        \,\middle|\,
        d_j(t)=1
    \right)
    &\le
    \exp\left(
        -\frac{(\mu_j^{k,t})^2}
        {8(H_{k,t}-1)}
    \right) \le
    \exp\left(
        -\frac{\underline\lambda^2}{8}H_{k,t}
    \right),
\end{align*}
where the last inequality follows from $\mu_{j}^{k,t} \geq \underline\lambda H_{k,t}$ and $\frac{H_{k,t}^2}{H_{k,t}-1} \geq H_{k,t}$.

\paragraph{Case 2: Wrong acceptance.}
Suppose the policy accepts but acceptance is not satisfying. Then $Q_j^{k,t}\ge\frac12\mu_j^{k,t}$ and $\widehat Q_j^{k,t}<1$, which imply
\begin{align}
    \sum_{\substack{i\in\mathcal J_k\\ i<j}}
    \left(
        \widehat A_i^{(k)}(t)
        -
        \mathbb E[\widehat A_i^{(k)}(t)]
    \right)
    >
    \frac{1}{2}\mu_j^{k,t}-1
    \ge
    \frac14\mu_j^{k,t}.
    \label{eq:wrong_acceptance_deviation}
\end{align}

Similarly, for $s=t,\ldots,\tau_k$, let $Y_s:=\sum_{i\in\mathcal J_k,\,i<j}d_i(s)$ and $p_{<j,s}:=\sum_{i\in\mathcal J_k,\,i<j}\lambda_{i,s}$. Conditional on $d_j(t)=1$, we have $Y_t=0$, so the current-period centered term is $-p_{<j,t}\le0$. Thus, \eqref{eq:wrong_acceptance_deviation} implies $\sum_{s=t+1}^{\tau_k}(Y_s-p_{<j,s})>\frac14\mu_j^{k,t}$. The future $Y_s$ remain independent Bernoulli random variables with means $p_{<j,s}$, so Hoeffding's inequality bounds the conditional probability of wrong acceptance by $\exp(-\underline\lambda^2 H_{k,t}/8)$ as well. Combining the two cases gives
\begin{align*}
    \mathbb P\left(
        \mathrm{Error}_{j,t}^{(k)}
        \,\middle|\,
        d_j(t)=1
    \right) \le
    2\exp\left(
        -\frac{\underline\lambda^2}{8}H_{k,t}
    \right).
\end{align*}
This proves \eqref{eq:direct_tracking_error_bound}.

Since at most one arrival occurs in period $t$, the tracking-error events
for different types are disjoint. Each requires an arrival of its type, so
\begin{align*}
    \mathbb P(\mathrm{Error}_t)
    &=\sum_{\substack{j\in\mathcal J_k\\\lambda_{j,t}>0}}
        \lambda_{j,t}\,
        \mathbb P\left(\mathrm{Error}_{j,t}^{(k)}\mid d_j(t)=1\right)\\
    &\le 2\exp\left(-\frac{\underline\lambda^2}{8}H_{k,t}\right)
        \sum_{j=1}^J\lambda_{j,t}
    \le 2\exp\left(-\frac{\underline\lambda^2}{8}H_{k,t}\right),
\end{align*}
where the first inequality uses \eqref{eq:direct_tracking_error_bound}
and the second uses $\sum_{j=1}^J\lambda_{j,t}\le1$. This proves
\eqref{eq:union_tracking_error_bound}.
\Halmos
\endproof

\section{Proofs for the Expiration Protection Policy}

In this section, we prove Lemmas~\ref{lem:heterogeneous_pointwise_count} and~\ref{lem:heterogeneous_EP_one_step}.

\subsection{Proof of Lemma~\ref{lem:heterogeneous_pointwise_count}}
\label{app:heterogeneous_pointwise_count}

For this and the following proofs, write $\vartheta_j:=\theta_j(e_j)$ and extend the definition of $\eta_j^{\tau_k}$ to periods $t\le e_j$ by
\begin{align*}
    \eta_j^t
    :=
    \vartheta_j
    +
    \sum_{s=t+1}^{e_j}
    \left(
        \sum_{i<j}\lambda_{i,s}
        +
        \frac12\lambda_{j,s}
    \right),
    \qquad t\in\{1,\ldots,e_j\}.
\end{align*}
For each phase $k$ and type $j\in\mathcal J_k$, define the rounded effective end-of-phase protection level, capped at $B$, by
\begin{align}
    g_j^{(k)}
    :=
    \min\left\{B,\left\lceil\eta_j^{\tau_k}\right\rceil\right\},
    \qquad k\in[K],\quad j\in\mathcal J_k.
    \label{eq:heterogeneous_integer_base}
\end{align}

\proof{Proof.}

Fix a phase $k$, an entering capacity $b\in\{0,\ldots,B\}$, and a cutoff $v\in[0,r_{\max}]$. Let $\pi$ denote the EP policy, given by Algorithm~\ref{Algo:GeneralPolicy} with $\zeta_j(t)=e_j$, run during phase $k$ from capacity $b$.
Write $B_t^{\mathrm{EP}}$ for the remaining capacity in this run, so $B_{\tau_{k-1}+1}^{\mathrm{EP}}=b$ and $B_{\tau_k+1}^{\mathrm{EP}}=b'$.

The proof has three steps. First, we identify exactly which policy decisions contribute to $N^*(v)-N^{\pi}(v)$. Second, we use the policy thresholds to derive deterministic bounds on the states at which these decisions can occur. Third, we derive the bound in \eqref{eq:heterogeneous_pointwise_separation} for each realization of the phase arrivals, then take expectations and apply Hoeffding's inequality to obtain \eqref{eq:heterogeneous_phasewise_binary_expected}.

\paragraph{Step 1: identify the decisions that contribute to the value-count difference.}

Call types with rewards above $v$ high and the remaining types low. Let $m$ denote the number of high types, which are $1,\ldots,m$ by the reward ordering, and let $h:=\sum_{u=1}^{b}\mathbf 1\{\nu_u^{(k)}>v\}$ count the continuation marginal values above $v$. For $\tau_{k-1}+1\le t\le\tau_k+1$, define
\begin{align*}
    H_t:=\sum_{s=t}^{\tau_k}\sum_{j\le m}d_j(s),
    \qquad
    Y_t:=B_t^{\mathrm{EP}}-H_t.
\end{align*}
Here, $H_t$ counts the high arrivals remaining in the phase, including period $t$. In particular, $Y_{\tau_{k-1}+1}=b-H_{\tau_{k-1}+1}$ and $Y_{\tau_k+1}=b'$. Rejecting a high arrival increases $Y_t$ by one, accepting a low arrival decreases it by one, and all other outcomes leave it unchanged.

There are $H_{\tau_{k-1}+1}+h$ arrival rewards and continuation marginal values above $v$ available to the benchmark, therefore 
\begin{align*}
    N^*(v)=\min\{b,H_{\tau_{k-1}+1}+h\} =H_{\tau_{k-1}+1}+\min\{Y_{\tau_{k-1}+1},h\}.
\end{align*}

The EP policy's accepted high arrivals are the total high arrivals minus the rejected ones, each corresponding to an upward step of $Y_t$. The $b'$ units of capacity carried forward contribute the continuation marginal values $\nu_1^{(k)},\ldots,\nu_{b'}^{(k)}$, of which $\min\{b',h\}$ exceed $v$. Using $Y_{\tau_k+1}=b'$, we obtain
\begin{align*}
    N^\pi(v) =H_{\tau_{k-1}+1}
      -\sum_{t\in\mathcal T_k}\mathbf 1\{Y_{t+1}>Y_t\}
      +\min\{Y_{\tau_k+1},h\}.
\end{align*}
Since $Y_t$ changes by at most one and $h$ is an integer,
\begin{align*}
    \min\{Y_{t+1},h\}-\min\{Y_t,h\} =
    \mathbf 1\{Y_{t+1}>Y_t,\ Y_t<h\}
    -\mathbf 1\{Y_{t+1}<Y_t,\ Y_t\le h\}.
\end{align*}
Subtracting the two counts and summing these increments yields
\begin{align}
    N^*(v)-N^\pi(v)
    =\sum_{t\in\mathcal T_k}
    \mathbf 1\{Y_{t+1}>Y_t,\ Y_t\ge h\} +\sum_{t\in\mathcal T_k}
    \mathbf 1\{Y_{t+1}<Y_t,\ Y_t\le h\}.
    \label{eq:heterogeneous_phasewise_wrong_step_identity}
\end{align}
Thus, only rejected high arrivals with $Y_t\ge h$ and accepted low arrivals with $Y_t\le h$ contribute to the count difference.

\paragraph{Step 2: use the policy thresholds to bound the states at which these decisions can occur.}

We now obtain an upper bound on $Y_t$ when a high arrival is rejected and a lower bound when a low arrival is accepted. Write $n_t:=\tau_k-t$ for the number of periods remaining after $t$, and set
\begin{align*}
    \bar g_j:=\min\{\lceil\eta_j^{\tau_k}\rceil,b\}
    =\min\{g_j^{(k)},b\},
    \qquad
    \omega_{j,t}:=\sum_{q<j}\lambda_{q,t}+\frac12\lambda_{j,t},
\end{align*}
so that $0\le\omega_{j,t}\le1$. Recall from \eqref{eq:heterogeneous_arrival_deviation} that $\Delta_v(t)$ is the number of high arrivals strictly after period $t$ through the end of phase $k$, minus its expectation.

By \eqref{eq:general_acceptance} and the definition of $\eta_j^{\tau_k}$, a type-$j$ arrival in period $t$ is accepted if and only if
\begin{align*}
    B_t^{\mathrm{EP}}
    \ge\eta_j^{\tau_k}
    +\sum_{s=t+1}^{\tau_k}\omega_{j,s}
    +\omega_{j,t}.
\end{align*}
On such an arrival, $\omega_{j,t}>0$, so the inequality also ensures $B_t^{\mathrm{EP}}\ge1$. Since the arriving type remains active through $\tau_k$, the active-interval assumption gives $\lambda_{j,s}\ge\underline\lambda$ for $s=t,\ldots,\tau_k$.

\textit{Rejected high arrival.}
Suppose that a type-$j$ arrival with $j\le m$ is rejected. By the acceptance rule above,
\begin{align*}
    B_t^{\mathrm{EP}}
    <\eta_j^{\tau_k}+\sum_{s=t+1}^{\tau_k}\omega_{j,s}+\omega_{j,t}.
\end{align*}
If $\eta_j^{\tau_k}<b$, then $\bar g_j=\lceil\eta_j^{\tau_k}\rceil\ge\eta_j^{\tau_k}$. Replacing $\eta_j^{\tau_k}$ by $\bar g_j$ and bounding the current-period demand adjustment $\omega_{j,t}$ by $1$ gives
\begin{align*}
    B_t^{\mathrm{EP}}
    <\bar g_j+\sum_{s=t+1}^{\tau_k}\omega_{j,s}+1.
\end{align*}
If $\eta_j^{\tau_k}\ge b$, then $\bar g_j=b$ and $B_t^{\mathrm{EP}}\le b$, so the same upper bound holds.

Since $H_t$ includes the current high arrival as well as the high arrivals in the remaining periods,
\begin{align*}
    H_t=1+\sum_{s=t+1}^{\tau_k}\sum_{q\le m}\lambda_{q,s}+\Delta_v(t).
\end{align*}
Subtracting $H_t$ from the upper bound on $B_t^{\mathrm{EP}}$ cancels the two $1$ terms. Using $Y_t=B_t^{\mathrm{EP}}-H_t$ and the definition of $\omega_{j,s}$, we obtain
\begin{align}
    Y_t
    &<\bar g_j
      -\sum_{s=t+1}^{\tau_k}
      \left(\frac12\lambda_{j,s}+\sum_{j<q\le m}\lambda_{q,s}\right)
      -\Delta_v(t)\notag\\
    &\le\bar g_j-\frac12\underline\lambda n_t-\Delta_v(t).
    \label{eq:heterogeneous_phasewise_high_barrier}
\end{align}
The last inequality uses $\lambda_{j,s}\ge\underline\lambda$ in each of the $n_t$ remaining periods.

\textit{Accepted low arrival.}
Suppose that a type-$j$ arrival with $j>m$ is accepted. The current arrival is not high, so
\begin{align*}
    H_t=\sum_{s=t+1}^{\tau_k}\sum_{q\le m}\lambda_{q,s}+\Delta_v(t).
\end{align*}
Subtracting this expression from the acceptance inequality gives
\begin{align}
    Y_t
    &\ge\eta_j^{\tau_k}+\omega_{j,t}
      +\sum_{s=t+1}^{\tau_k}
      \left(\sum_{m<q<j}\lambda_{q,s}+\frac12\lambda_{j,s}\right)
      -\Delta_v(t)\notag\\
    &>\bar g_j-1+\frac12\underline\lambda n_t-\Delta_v(t),
    \label{eq:heterogeneous_phasewise_low_barrier}
\end{align}
where the last inequality uses $\eta_j^{\tau_k}>\lceil\eta_j^{\tau_k}\rceil-1\ge\bar g_j-1$, $\omega_{j,t}\ge0$, and $\lambda_{j,s}\ge\underline\lambda$ for the remaining periods.

\paragraph{Step 3: derive the bound for each realization and take expectations.}

Recall that $h=\sum_{u=1}^{b}\mathbf 1\{\nu_u^{(k)}>v\}$ counts the continuation marginal values above $v$ among $\nu_1^{(k)},\ldots,\nu_b^{(k)}$. By \eqref{eq:heterogeneous_d_definition}, for each $j\in\mathcal J_k$, we have $\bar g_j-h\le\xi(v)$ if $j\le m$ and $h-\bar g_j\le\xi(v)$ if $j>m$.  If a rejected high arrival of type $j$ contributes to \eqref{eq:heterogeneous_phasewise_wrong_step_identity}, then $Y_t\ge h$, so \eqref{eq:heterogeneous_phasewise_high_barrier} gives
\begin{align*}
    \frac12\underline\lambda n_t
    <\bar g_j-h-\Delta_v(t)
    \le\xi(v)+|\Delta_v(t)|.
\end{align*}
If an accepted low arrival of type $j$ contributes, then $Y_t\le h$, and \eqref{eq:heterogeneous_phasewise_low_barrier} gives
\begin{align*}
    \frac12\underline\lambda n_t
    <h-\bar g_j+1+\Delta_v(t)
    \le\xi(v)+1+|\Delta_v(t)|.
\end{align*}
Thus, every decision counted in \eqref{eq:heterogeneous_phasewise_wrong_step_identity} satisfies
\begin{equation}
    \frac12\underline\lambda n_t
    <\xi(v)+1+|\Delta_v(t)|.
    \label{eq:heterogeneous_count_contribution_condition}
\end{equation}

For a contributing decision with $|\Delta_v(t)|\le\underline\lambda n_t/8$, \eqref{eq:heterogeneous_count_contribution_condition} implies
\begin{align*}
    \frac38\underline\lambda n_t<\xi(v)+1,
    \qquad\text{hence}\qquad
    n_t<\frac{8(\xi(v)+1)}{3\underline\lambda}.
\end{align*}
Since $n_t=\tau_k-t$ takes each integer from $0$ to $L_k-1$ once, at most $\lceil 8(\xi(v)+1)/(3\underline\lambda)\rceil$ periods satisfy this last bound. Every other period contributing to the count difference must satisfy $|\Delta_v(t)|>\underline\lambda n_t/8$. There is at most one contributing decision per period, so \eqref{eq:heterogeneous_phasewise_wrong_step_identity} gives \eqref{eq:heterogeneous_pointwise_separation}.

For $n_t\ge1$, $\Delta_v(t)$ is the sum of $n_t$ independent high-arrival indicators minus their expectations. Hoeffding's inequality therefore gives
\begin{equation}
    \mathbb P\!\left(|\Delta_v(t)|>\frac18\underline\lambda n_t\right)
    \le 2\exp\!\left(-\frac{\underline\lambda^2}{32}n_t\right).
    \label{eq:heterogeneous_count_deviation_probability}
\end{equation}
For $n_t=0$, we have $\Delta_v(t)=0$, so this event cannot occur. Summing over periods yields
\begin{align*}
    \sum_{t\in\mathcal T_k}
    \mathbb P\!\left(|\Delta_v(t)|>\frac18\underline\lambda n_t\right)
    \le 2\sum_{n=1}^{\infty}\exp\!\left(-\frac{\underline\lambda^2}{32}n\right)
    =\frac{2}{\exp(\underline\lambda^2/32)-1}.
\end{align*}
With $b$ fixed, taking expectations over the phase arrivals in \eqref{eq:heterogeneous_pointwise_separation} and using $\lceil x\rceil\le x+1$ proves \eqref{eq:heterogeneous_phasewise_binary_expected}.
\Halmos
\endproof

\subsection{Proof of Lemma~\ref{lem:heterogeneous_EP_one_step}}
\label{app:heterogeneous_EP_one_step}

\proof{Proof of Lemma~\ref{lem:heterogeneous_EP_one_step}.}

Recall from \eqref{eq:heterogeneous_integer_base} that $g_j^{(k)}=\min\{B,\lceil\eta_j^{\tau_k}\rceil\}$ is the effective end-of-phase protection level rounded up to an integer and capped at $B$.

Fix a phase $k$, a type $j\in\mathcal J_k$ with $k<\kappa_j$, and a phase $\ell\in\{k+1,\ldots,\kappa_j\}$. Define the increment
\begin{align*}
    \sigma:=\left\lceil\eta_j^{\tau_{\ell-1}}\right\rceil
    -\left\lceil\eta_j^{\tau_\ell}\right\rceil.
\end{align*}
The increment $\sigma$ is a nonnegative integer computed before capping the protection levels at $B$, so $g_j^{(\ell-1)}=\min\{B,g_j^{(\ell)}+\sigma\}$.

By the definition of the effective end-of-phase protection levels,
\begin{align*}
    \eta_j^{\tau_{\ell-1}}-\eta_j^{\tau_\ell}
    =\mathbb E\left[\sum_{i<j}A_i^{(\ell)}\right]+\frac12\mathbb E[A_j^{(\ell)}].
\end{align*}
Since $|\lceil x+y\rceil-\lceil x\rceil-y|<1$, we obtain
\begin{equation}
    \left|\sigma-\left(\mathbb E\left[\sum_{i<j}A_i^{(\ell)}\right]
    +\frac12\mathbb E[A_j^{(\ell)}]\right)\right|<1.
    \label{eq:heterogeneous_increment_condition}
\end{equation}

For each realization $\boldsymbol A^{(\ell)}$ and capacity $b\in\{0,\ldots,B\}$, define
\begin{align*}
    G_{\boldsymbol A^{(\ell)}}(b)
    :=\max_{\tilde b\in\{0,\ldots,b\}}
    \left\{P(b-\tilde b,\boldsymbol A^{(\ell)})
    +\Psi_{\ell+1}(\tilde b)\right\}-r_jb.
\end{align*}
Here, $G_{\boldsymbol A^{(\ell)}}(b)$ measures the net value of retaining $b$ units for phase $\ell$ and later phases, using the type-$j$ reward $r_j$ as the opportunity cost per unit of capacity.
The benchmark recursion \eqref{eq:phase_aware_recursion} gives $\mathbb E[G_{\boldsymbol A^{(\ell)}}(b)]=\Psi_\ell(b)-r_jb$. Hence, by \eqref{eq:heterogeneous_protection_level_value_loss}, the protection-level value loss at $g_j^{(\ell-1)}$ satisfies
\begin{align*}
    \mathcal L_{\Psi_\ell}^j(g_j^{(\ell-1)})
    =\max_{b\in\{0,\ldots,B\}}
    \mathbb E\!\left[G_{\boldsymbol A^{(\ell)}}(b)
    -G_{\boldsymbol A^{(\ell)}}(g_j^{(\ell-1)})\right].
\end{align*}
This loss is the reduction in expected net value from using $g_j^{(\ell-1)}$ rather than the benchmark's phase-specific protection level $\theta_j(\tau_{\ell-1})$.

We proceed in four steps. First, we identify arrival realizations for which $\sigma$ maximizes $P(q,\boldsymbol A^{(\ell)})-r_jq$. Second, for these realizations, we bound the difference between the two $G$ terms by $\mathcal L_{\Psi_{\ell+1}}^j(g_j^{(\ell)})$. Third, we bound the additional contribution from the remaining realizations and take expectations. Finally, we use Hoeffding's inequality to bound the probability of this complementary event.

\paragraph{Step 1: identify when $\sigma$ is optimal for the realized current-phase problem.}

Define the event
\begin{align*}
    \mathcal G_\ell^j
    :=\left\{\sum_{i<j}A_i^{(\ell)}\le\sigma\le \sum_{i\le j}A_i^{(\ell)}\right\}.
\end{align*}
We use \eqref{eq:offline_value_def} to define $P(q,\boldsymbol A^{(\ell)})$ for any nonnegative integer $q$, including $q>B$. On $\mathcal G_\ell^j$, the first $\sigma$ arrival rewards in nonincreasing order are at least $r_j$, and all remaining rewards are at most $r_j$. Since $P$ is constant beyond the total realized demand, $P(q,\boldsymbol A^{(\ell)})-r_jq$ is nondecreasing up to $\sigma$ and nonincreasing thereafter. Thus,
\begin{equation}
    P(q,\boldsymbol A^{(\ell)})-r_jq
    \le P(\sigma,\boldsymbol A^{(\ell)})-r_j\sigma,
    \qquad q\in\mathbb Z_{\ge0}.
    \label{eq:heterogeneous_good_event_marginal}
\end{equation}

\paragraph{Step 2: bound the difference between the two $G$ terms on $\mathcal G_\ell^j$.}

Fix a realization in $\mathcal G_\ell^j$ and a capacity $b\in\{0,\ldots,B\}$. Consider any nonnegative integers $q,y$ with $q+y=b$, allocating $q$ units to phase $\ell$ and carrying $y$ units forward. By the definition of the protection-level value loss,
\begin{align*}
    \Psi_{\ell+1}(y)-r_jy
    \le \Psi_{\ell+1}(g_j^{(\ell)})-r_jg_j^{(\ell)}
    +\mathcal L_{\Psi_{\ell+1}}^j(g_j^{(\ell)}).
\end{align*}

\textbf{Case 1: $g_j^{(\ell)}+\sigma\le B$.}
Here $g_j^{(\ell-1)}=g_j^{(\ell)}+\sigma$, so allocating $\sigma$ units to phase $\ell$ and carrying $g_j^{(\ell)}$ units forward is feasible at capacity $g_j^{(\ell-1)}$. Combining \eqref{eq:heterogeneous_good_event_marginal} with the preceding inequality gives
\begin{align*}
    P(q,\boldsymbol A^{(\ell)})+\Psi_{\ell+1}(y)-r_jb 
    & \le P(\sigma,\boldsymbol A^{(\ell)})-r_j\sigma
    +\Psi_{\ell+1}(g_j^{(\ell)})-r_jg_j^{(\ell)}
    +\mathcal L_{\Psi_{\ell+1}}^j(g_j^{(\ell)})
    \\
    & \le G_{\boldsymbol A^{(\ell)}}(g_j^{(\ell-1)})
    +\mathcal L_{\Psi_{\ell+1}}^j(g_j^{(\ell)}),
\end{align*}
where the second inequality follows from the definition of $G_{\boldsymbol A^{(\ell)}}$.

\textbf{Case 2: $g_j^{(\ell)}+\sigma>B$.}
Here $g_j^{(\ell-1)}=B$. If $B-y\le\sigma$, then $q=b-y\le B-y\le\sigma$. The monotonicity established in Step~1 therefore gives
\begin{align*}
    P(q,\boldsymbol A^{(\ell)})+\Psi_{\ell+1}(y)-r_jb
    &\le P(B-y,\boldsymbol A^{(\ell)})+\Psi_{\ell+1}(y)-r_jB\\
    &\le G_{\boldsymbol A^{(\ell)}}(B).
\end{align*}
If $B-y>\sigma$, then $0\le y<B-\sigma<g_j^{(\ell)}\le B$. By Lemma~\ref{lem:psi_discrete_concavity}, the function $u\mapsto\Psi_{\ell+1}(u)-r_ju$ is discretely concave, so
\begin{align*}
    \Psi_{\ell+1}(y)-r_jy
    -\bigl[\Psi_{\ell+1}(B-\sigma)-r_j(B-\sigma)\bigr]
    & \le\bigl[\Psi_{\ell+1}(y)-r_jy
    -\Psi_{\ell+1}(g_j^{(\ell)})+r_jg_j^{(\ell)}\bigr]^+
    \\
    & \le\mathcal L_{\Psi_{\ell+1}}^j(g_j^{(\ell)}).
\end{align*}
Together with \eqref{eq:heterogeneous_good_event_marginal}, this yields
\begin{align*}
    P(q,\boldsymbol A^{(\ell)})+\Psi_{\ell+1}(y)-r_jb
    & \le P(\sigma,\boldsymbol A^{(\ell)})
    +\Psi_{\ell+1}(B-\sigma)-r_jB
    +\mathcal L_{\Psi_{\ell+1}}^j(g_j^{(\ell)})\\
    & \le G_{\boldsymbol A^{(\ell)}}(B)
    +\mathcal L_{\Psi_{\ell+1}}^j(g_j^{(\ell)}),
\end{align*}
since allocating $\sigma$ units to phase $\ell$ and carrying $B-\sigma$ units forward is feasible at capacity $B$.

In both cases, maximizing over the feasible allocations $q+y=b$ gives
\begin{align*}
    G_{\boldsymbol A^{(\ell)}}(b)
    -G_{\boldsymbol A^{(\ell)}}(g_j^{(\ell-1)})
    \le\mathcal L_{\Psi_{\ell+1}}^j(g_j^{(\ell)}),
    \qquad b\in\{0,\ldots,B\}.
\end{align*}

\paragraph{Step 3: bound the additional expected loss from the complementary event.}

We now bound the difference between the two $G$ terms on $(\mathcal G_\ell^j)^c$. Since at most one arrival occurs per period, the expected demand adjustment in \eqref{eq:heterogeneous_increment_condition} is at most $L_\ell$. That inequality gives $\sigma<L_\ell+1$; since $\sigma$ is an integer,
\begin{align*}
    0\le g_j^{(\ell-1)}-g_j^{(\ell)}\le\sigma\le L_\ell.
\end{align*}
Fix a realization in $(\mathcal G_\ell^j)^c$ and a capacity $b\in\{0,\ldots,B\}$. Choose an optimal allocation $q+y=b$ in $G_{\boldsymbol A^{(\ell)}}(b)$, with $q$ units allocated to phase $\ell$ and $y$ carried forward. By monotonicity of $\Psi_{\ell+1}$, any capacity allocated beyond the total realized demand can instead be carried forward without reducing the objective. We may therefore choose
\begin{align*}
    0\le q\le\sum_{i=1}^J A_i^{(\ell)}\le L_\ell.
\end{align*}
At capacity $g_j^{(\ell-1)}$, it is feasible to allocate $g_j^{(\ell-1)}-g_j^{(\ell)}$ units to phase $\ell$ and carry $g_j^{(\ell)}$ forward. Comparing with this allocation gives
\begin{align*}
    G_{\boldsymbol A^{(\ell)}}(b)
    -G_{\boldsymbol A^{(\ell)}}(g_j^{(\ell-1)})
    & \le\bigl[P(q,\boldsymbol A^{(\ell)})-r_jq\bigr] -\bigl[P(g_j^{(\ell-1)}-g_j^{(\ell)},\boldsymbol A^{(\ell)})
    -r_j(g_j^{(\ell-1)}-g_j^{(\ell)})\bigr]\\
    & \quad+\bigl[\Psi_{\ell+1}(y)-r_jy\bigr]
    -\bigl[\Psi_{\ell+1}(g_j^{(\ell)})-r_jg_j^{(\ell)}\bigr]\\
    & \le 2r_{\max}L_\ell
    +\mathcal L_{\Psi_{\ell+1}}^j(g_j^{(\ell)}),
\end{align*}
where the first inequality follows from the definition of $G_{\boldsymbol A^{(\ell)}}$, using an optimal allocation at capacity $b$ and a feasible allocation at capacity $g_j^{(\ell-1)}$.
The second follows because the first two bracketed terms each lie in $[-r_{\max}L_\ell,r_{\max}L_\ell]$, as both current-phase allocations are at most $L_\ell$, while the difference between the last two terms is at most $\mathcal L_{\Psi_{\ell+1}}^j(g_j^{(\ell)})$ by definition.

Combining this bound with Step~2 gives, for every realization and every $b\in\{0,\ldots,B\}$,
\begin{align*}
    G_{\boldsymbol A^{(\ell)}}(b)
    -G_{\boldsymbol A^{(\ell)}}(g_j^{(\ell-1)})
    \le\mathcal L_{\Psi_{\ell+1}}^j(g_j^{(\ell)})
    +2r_{\max}L_\ell\mathbf 1\{(\mathcal G_\ell^j)^c\}.
\end{align*}
Taking expectations and then maximizing over $b$, using the representation of the loss at the start of the proof, yields
\begin{align}
    \mathcal L_{\Psi_\ell}^j(g_j^{(\ell-1)})
    \le\mathcal L_{\Psi_{\ell+1}}^j(g_j^{(\ell)})
    +2r_{\max}L_\ell\mathbb P((\mathcal G_\ell^j)^c).
    \label{eq:heterogeneous_loss_before_concentration}
\end{align}

\paragraph{Step 4: bound the probability of the complementary event.}

By the definition of $\mathcal G_\ell^j$,
\begin{align*}
    (\mathcal G_\ell^j)^c
    =\left\{\sum_{i<j}A_i^{(\ell)}>\sigma\right\}
    \cup\left\{\sum_{i\le j}A_i^{(\ell)}<\sigma\right\}.
\end{align*}
Let $\mu:=\mathbb E[A_j^{(\ell)}]$. On the first and second events, respectively, the integrality of the arrival counts and $\sigma$, together with \eqref{eq:heterogeneous_increment_condition}, gives
\begin{align*}
    \sum_{i<j}A_i^{(\ell)}
    &\ge\sigma+1
    >\mathbb E\left[\sum_{i<j}A_i^{(\ell)}\right]+\frac12\mu,\\
    \sum_{i\le j}A_i^{(\ell)}
    &\le\sigma-1
    <\mathbb E\left[\sum_{i\le j}A_i^{(\ell)}\right]-\frac12\mu.
\end{align*}
Thus, each event requires a deviation of more than $\mu/2$ from the corresponding expected arrival count. Each count is a sum of $L_\ell$ arrival indicators that are independent across periods. Hoeffding's inequality bounds the probability of each event by $\exp[-\mu^2/(2L_\ell)]$.

Since $j\in\mathcal J_k$ and $k<\ell\le\kappa_j$, type $j$ remains active throughout phase $\ell$. By the active-interval assumption,
\begin{align*}
    \mu=\sum_{t\in\mathcal T_\ell}\lambda_{j,t}
    \ge\underline\lambda L_\ell.
\end{align*}
A union bound therefore gives
\begin{align*}
    \mathbb P((\mathcal G_\ell^j)^c)
    \le 2\exp\left(-\frac{\mu^2}{2L_\ell}\right)
    \le 2\exp\left(-\frac{\underline\lambda^2}{2}L_\ell\right).
\end{align*}
Substituting this bound into \eqref{eq:heterogeneous_loss_before_concentration} yields
\begin{align*}
    \mathcal L_{\Psi_\ell}^j(g_j^{(\ell-1)})
    \le\mathcal L_{\Psi_{\ell+1}}^j(g_j^{(\ell)})
    +4r_{\max}L_\ell
    \exp\left(-\frac{\underline\lambda^2}{2}L_\ell\right),
\end{align*}
which proves the lemma by the definition of $g_j^{(\ell)}$ in \eqref{eq:heterogeneous_integer_base}.

\Halmos
\endproof

\section{Effect of Approximating Arrival Probabilities}
\label{app:arrival_approximation}

In this section, we bound the change in expected reward when a protection policy constructed using modified arrival probabilities is applied to the original arrival process.
Fix the horizon $T$, initial capacity $B$, and rewards $r_1,\ldots,r_J$.
Consider two arrival models with probabilities $\lambda_{j,t}$ and $\tilde\lambda_{j,t}$, respectively, with nonnegative probabilities summing to at most one in each period and independent arrival outcomes across periods.
Write $\mathbb E_\lambda$ and $\mathbb E_{\tilde\lambda}$ for expectations under the two models.

Suppose the modified model ($\tilde\lambda$) satisfies the active-interval assumption, and let $\tilde\pi$ be either the PP policy or the EP policy constructed using $\tilde\lambda$.
All protection reference times, protection levels, and demand adjustments are computed under this modified model and remain fixed when the policy is applied to either arrival process.
For $\tilde\lambda_{j,t}>0$, let $h_{j,t}$ denote the acceptance threshold in \eqref{eq:general_acceptance}, i.e.,
\begin{align*}
    h_{j,t}
    :=
    \theta_j(\zeta_j(t))
    +\sum_{s=t}^{\zeta_j(t)}
    \left(
        \sum_{i<j}\tilde\lambda_{i,s}
        +\frac12\tilde\lambda_{j,s}
    \right).
\end{align*}
Under the original model ($\lambda$), a type-$j$ arrival may occur in a period $t$ even when $\tilde\lambda_{j,t}=0$.
We extend $\tilde\pi$ to such arrivals by rejecting them.
Setting $h_{j,t}=+\infty$ in these periods ensures that no finite remaining capacity satisfies $B_t\ge h_{j,t}$.
Thus, under either arrival model, $\tilde\pi$ accepts a type-$j$ arrival if and only if $B_t\ge1$ and $B_t\ge h_{j,t}$.

The following proposition bounds the difference between the expected rewards of this fixed policy under the two arrival models by the total change in arrival probabilities.

\begin{proposition}[Expected reward under arrival approximation]
\label{prop:arrival_approximation_reward}
The policy $\tilde\pi$ satisfies
\begin{align}
    \left|
        \mathbb E_\lambda[\mathrm{Reward}(\tilde\pi)]
        -\mathbb E_{\tilde\lambda}[\mathrm{Reward}(\tilde\pi)]
    \right|
    \le
    r_{\max}\sum_{t=1}^T\sum_{j=1}^J
        |\lambda_{j,t}-\tilde\lambda_{j,t}|
    =r_{\max}\varepsilon.
    \label{eq:arrival_approximation_reward_bound}
\end{align}
\end{proposition}

\proof{Proof.}
First, let $Y_t$ and $\tilde Y_t$ denote the period-$t$ arrival outcomes under $\lambda$ and $\tilde\lambda$, respectively.
Both take values in $\{0,1,\ldots,J\}$, where $0$ denotes no arrival and $j\ge1$ denotes a type-$j$ arrival.
Write $\lambda_{0,t}:=1-\sum_{j=1}^J\lambda_{j,t}$ and $\tilde\lambda_{0,t}:=1-\sum_{j=1}^J\tilde\lambda_{j,t}$ for the no-arrival probabilities.
Couple the two outcomes with their prescribed distributions so that
\begin{align*}
    \mathbb P(Y_t=\tilde Y_t=j)
    =\min\{\lambda_{j,t},\tilde\lambda_{j,t}\},
    \qquad j=0,\ldots,J.
\end{align*}
Choose these pairs independently across periods. Then
\begin{align*}
    \mathbb P(Y_t\ne\tilde Y_t)
    &=1-\sum_{j=0}^J\min\{\lambda_{j,t},\tilde\lambda_{j,t}\}\\
    &=\frac12\sum_{j=0}^J|\lambda_{j,t}-\tilde\lambda_{j,t}|\\
    &=\frac12\left(
        \sum_{j=1}^J|\lambda_{j,t}-\tilde\lambda_{j,t}|
        +\left|\sum_{j=1}^J(\lambda_{j,t}-\tilde\lambda_{j,t})\right|
    \right)\\
    &\le \sum_{j=1}^J|\lambda_{j,t}-\tilde\lambda_{j,t}|.
\end{align*}
Here, the second equality uses $2\min\{a,b\}=a+b-|a-b|$, the third substitutes the no-arrival probabilities, and the inequality follows from the triangle inequality.

Second, to bound the reward difference caused by these mismatches, we compare the policy's decisions on arrival sequences that differ in only one period.
For a deterministic sequence $\boldsymbol y=(y_1,\ldots,y_T)\in\{0,1,\ldots,J\}^T$, let $R(\boldsymbol y)$ denote the total reward earned by $\tilde\pi$, and let $B_t(\boldsymbol y)$ denote its remaining capacity at the beginning of period $t$.
Consider two such sequences $\boldsymbol y$ and $\boldsymbol z$ that differ only in period $u$.
The policy makes the same decisions before period $u$, so $B_u(\boldsymbol y)=B_u(\boldsymbol z)$.
There are three cases.
\begin{enumerate}
    \item The policy accepts nothing in period $u$ under either sequence.
    The capacities remain equal, and all subsequent arrivals and decisions coincide. Hence, $R(\boldsymbol y)=R(\boldsymbol z)$.

    \item Both arrivals are accepted in period $u$.
    The capacities are again equal after that period, and all subsequent decisions coincide. Thus,
    \begin{align*}
        |R(\boldsymbol y)-R(\boldsymbol z)|
        =|r_{y_u}-r_{z_u}|\le r_{\max}.
    \end{align*}

    \item Exactly one arrival is accepted in period $u$.
    Without loss of generality, suppose it is the arrival under $\boldsymbol y$, earning reward $r_{y_u}$. Then
    \begin{align*}
        B_{u+1}(\boldsymbol z)-B_{u+1}(\boldsymbol y)=1.
    \end{align*}
    From period $u+1$ onward, arrivals and thresholds are identical under the two sequences.
    Whenever the policy accepts with less capacity, it also accepts with more capacity, so the capacities remain one unit apart until the first period $v>u$, if any, in which the decisions differ.
    In that period, the policy accepts only under $\boldsymbol z$, earning reward $r_{z_v}$, and
    \begin{align*}
        B_{v+1}(\boldsymbol z)=B_{v+1}(\boldsymbol y).
    \end{align*}
    All subsequent decisions then coincide. Therefore,
    \begin{align*}
        R(\boldsymbol y)-R(\boldsymbol z)
        =
        \begin{cases}
            r_{y_u}-r_{z_v}, & \text{if such a period $v$ exists},\\
            r_{y_u}, & \text{otherwise}.
        \end{cases}
    \end{align*}
    In either case, $|R(\boldsymbol y)-R(\boldsymbol z)|\le r_{\max}$.
\end{enumerate}
Thus, changing one arrival outcome changes total reward by at most $r_{\max}$.

Finally, we extend this one-period comparison to arbitrary arrival sequences $\boldsymbol y$ and $\boldsymbol z$.
Starting from $\boldsymbol y$, replace $y_1$ by $z_1$, then $y_2$ by $z_2$, and continue through period $T$.
The sequence after the first $t$ replacements is
\begin{align*}
    \boldsymbol y^{(t)}
    :=(z_1,\ldots,z_t,y_{t+1},\ldots,y_T),
    \qquad t=0,\ldots,T.
\end{align*}
Thus, $\boldsymbol y^{(0)}=\boldsymbol y$ and $\boldsymbol y^{(T)}=\boldsymbol z$.
The sequences $\boldsymbol y^{(t-1)}$ and $\boldsymbol y^{(t)}$ agree in every period except possibly period $t$, where $y_t$ is replaced by $z_t$.
Each replacement therefore changes total reward by at most $r_{\max}$, and changes nothing if $y_t=z_t$.
Summing these bounds and applying the triangle inequality gives
\begin{align*}
    |R(\boldsymbol y)-R(\boldsymbol z)|
    &\le \sum_{t=1}^T
        |R(\boldsymbol y^{(t-1)})-R(\boldsymbol y^{(t)})|\\
    &\le r_{\max}\sum_{t=1}^T\mathbf1\{y_t\ne z_t\}.
\end{align*}
Applying this bound to the coupled sequences $\boldsymbol Y=(Y_1,\ldots,Y_T)$ and $\tilde{\boldsymbol Y}=(\tilde Y_1,\ldots,\tilde Y_T)$ and taking expectations yields
\begin{align*}
    \left|
        \mathbb E_\lambda[\mathrm{Reward}(\tilde\pi)]
        -\mathbb E_{\tilde\lambda}[\mathrm{Reward}(\tilde\pi)]
    \right|
    &\le \mathbb E\bigl[|R(\boldsymbol Y)-R(\tilde{\boldsymbol Y})|\bigr]\\
    &\le r_{\max}\sum_{t=1}^T\mathbb P(Y_t\ne\tilde Y_t)\\
    &\le r_{\max}\varepsilon.
\end{align*}
\Halmos
\endproof

\end{document}